\documentclass[10pt]{article}

\usepackage[margin=1in]{geometry}
\usepackage{graphicx,enumitem}
\usepackage[usenames,dvipsnames]{xcolor}
\usepackage{amsmath}
\usepackage{amssymb}
\usepackage{amsthm}
\usepackage{multirow}
\usepackage[multiple]{footmisc}
\usepackage{cancel}
\usepackage{setspace}
\usepackage{comment}
\usepackage{caption}
\usepackage{mathrsfs}
\usepackage{tikz}
\DeclareMathOperator*{\argmax}{arg\,max}
\DeclareMathOperator*{\argmin}{arg\,min}

\newcommand{\PP}{{\mathbb P}}
\newcommand{\EE}{{\mathbb E}}
\newcommand{\FF}{{\mathcal F}}

\newcommand\hq{{\hat{q}}}

\newcommand\up{{\underline{p}}}
\newcommand\op{{\bar{p}}}

\newcommand{\defref}[1]{Definition~\ref{#1}}
\newcommand{\thmref}[1]{Theorem~\ref{#1}}
\newcommand{\proref}[1]{Proposition~\ref{#1}}

\newcommand{\lemref}[1]{Lemma~\ref{#1}}

\def\E{\mathbb{E}}

\numberwithin{equation}{section}
\newtheorem{theorem}{Theorem}[section]
\newtheorem{lemma}[theorem]{Lemma}
\newtheorem{proposition}[theorem]{Proposition}

\newtheorem{definition}[theorem]{Definition}

\newtheorem{remark}[theorem]{Remark}

\usepackage[toc,title,page]{appendix}

\begin{document}
\onehalfspacing
\title{A reputation game on cyber-security and cyber-risk calibration}
	
\author{Kookyoung Han\footnote{Ulsan National Institute of Science and Technology, Email: kyhan@unist.ac.kr} \, and Jin Hyuk Choi\footnote{Ulsan National Institute of Science and Technology, Email: jchoi@unist.ac.kr}}
\date{\today}
	
	
\maketitle
	
	\begin{abstract}
	To analyze strategic interactions arising in the cyber-security context, we develop a new reputation game model in which an attacker can pretend to be a normal user and a defender may have to announce attack detection at a certain point of time without knowing whether he has been attacked. We show the existence and uniqueness of sequential equilibrium in Markov strategies, and  explicitly characterize the players' equilibrium strategies. Using our model, we suggest empirical and theoretical ways of calibrating the attack probability, which is an important element of cyber-risks.
	\end{abstract}
	




\section{Introduction}
Modern technologies have been relying more and more on networks such as the Internet of Things, mobile networks, and cyber-physical systems. At the same time, many organizations and institutions have been suffering from bigger threats of cyber-attacks such as advanced persistent threats (APTs). Recent examples of APTs are the SolarWinds hack in 2020, Microsoft Exchange Server data breach in 2021, and Double Dragon (APT41). The first two examples are known to operate over a few months, and the last example is known to operation over a few years.

	There are several characteristics of APTs. First, it is difficult to prevent APTs because APT actors utilize various tools such as zero-day attacks, unfixed vulnerabilities of a system, and even social engineering. Second, it is not easy to detect APTs. APT actors steal a small amount of data, pretending to be normal users. Third, APTs are carried out for a long period of time. Since it is difficult to prevent and detect APTs, it takes a long time to become aware of APTs. Fourth, APT actors adjust their activities based on circumstances. These characteristics of APTs have two implications. First, it is difficult to discern malicious hackers from innocent users. Especially, this difficulty is persistent over time. Second, false alarms have to be involved in detecting APTs since it is too late if a security manager waits until she collects hard evidences and fully realizes cyber-attacks as a fact.

	APTs consist of multiple stages. At the reconnaissance stage, ATP actors lure users or employees of a targeted company. Once some are lured and infected, APT actors use the infected hosts as a foothold, and escalate their privileges to obtain an access to servers of the target company. Then, APT actors continuously and slowly steal data from the servers. A security manager estimates the likelihood of being attacked using a certain countermeasure. The security manager take actions such as shutting down the servers once she is sufficiently suspicious of cyber-attacks.

Motivated by the characteristics of APTs, we consider a dynamic game in which a defender tries to detect cyber-attacks, suffering from persistent private information and false alarm costs.  Persistent private information has been studied in reputation games and asset pricing models. Firms can pretend to be a commitment type to threaten potential entrants \cite{KrepsWilson1982, MilgromRoberts1982}. In \cite{Cripps2004}, it is shown that in a large class of repeated games, reputation effect eventually disappears. In \cite{Faingold}, reputation game in continuous time setting is studied when a group of small players faces a large commitment type player. Regarding asset pricing models, \cite{AndersonSmith2013, BackBaruch, Caldentey, campi2007insider, ccetin2018financial, CHOI201922, CollinFos2016, Kyle1985} investigate how equilibrium asset price dynamics is derived by informed trader's trading strategies.{\footnote{
Asset pricing has also been studied in models without  private information. For instance, \cite{choi2015taylor, weston2018existence, WesZit18, Zit12} study dynamic asset pricing in Radner equilibrium. 
} 

The defender in our model plays against a suspect who can be either an attacker who dynamically chooses actions or an innocent user who repeats the same action over time. Observing noisy signals of the suspect's actions, the defender can decide to see whether the suspect is the attacker or innocent by inspecting the suspect. The defender expects to incur potential damages due to cyber-attacks before inspection, but he incurs the false alarm cost when the suspect turns out to be innocent after inspection. This aspect is a difference between our model and the studies mentioned in the previous paragraph because there is no way that players in the studies reveal private information. However, the defender in our model can reveal private information although he is penalized for false detection.

We explicitly solve for sequential equilibrium in Markov strategies with the suspicion level, which is the posterior probability that the suspect is the attacker based on noisy observations of the suspect's actions. 
Our analysis shows that the attacker's actions are weakened as the suspicion level increases and that the defender begins inspection in equilibrium only if the suspicion level exceeds a certain threshold. In addition to characterization of equilibrium, we propose an empirical way of estimating the initial probability of the suspect being the attacker based on data that only indicate whether a user is inspected or not. We also propose a theoretical estimation on the initial probability of cyber-attacks, assuming that the attacker can choose the attack probability right before the game. Two methods that we propose would be useful because actual estimation may be neither available nor reliable due to lack of actual data.

The rest of this paper is organized as follows. In Section 2, we formally describe the model. In Section 3, we explicitly characterize sequential equilibrium of the model and provide comparative statics of equilibrium strategies. The proof of the main theorem is provided in Section 4.  In Section 5, we propose empirical and theoretical estimations on the initial probability of cyber-attacks. In Section 6, we calibrate the model parameters using a report on data breach. 
Section 7 illustrates graphical results and numerical simulations. Section 8 summarizes this paper and suggests several extensions of our model for future research.

\section{The Model}

We consider a continuous-time game between two risk-neutral players, a suspect and a defender. The suspect's type is a random variable $\theta$ taking value in $\{0,1\}$. The suspect is an attacker ($\theta = 1$) with probability $q_0 \in (0,1)$ or innocent ($\theta = 0$) with the complementary probability $1-q_0$. The suspect knows the true value of $\theta$, whereas the defender does not. 

The attacker chooses attack intensity $0 \leq \Delta_t \leq M$ at every moment $t \geq 0$ in time, where the constant $M>0$ is the upper bound of attack intensities. The innocent type always chooses zero attack intensity. One interpretation of attack intensities is the amount of data that the attacker steals at every moment in time. The defender chooses whether to block the suspect or not at every moment in time. Once the defender blocks the suspect, the game ends and the true value of $\theta$ is publicly revealed.\footnote{This assumption can be thought of as a circumstance in which after blocking a suspected user, a defender begins a thorough inspection that results in hard evidence on the identity of the suspected user.}

The defender does not directly observe the suspect's attack intensity. Instead, the defender observes the signal process $(Y_t)_{t \geq 0}$, which is noisy observations of the suspect's attack intensities. We assume that the signal process obeys the following stochastic differential equation (SDE):
\begin{equation}
\begin{split}\label{Y}
dY_t = \Delta_t {1}_{\{\theta=1\}} dt + \sigma \, dW_t,
\end{split}
\end{equation}
where $1_{\{\theta=1\}}$ is the indicator function, $\sigma$ is a strictly positive constant, and $(W_t)_{t\geq0}$ is a standard Brownian motion independent of $\theta$. The signal process $(Y_t)_{t \geq 0}$ is public information. That is, the players observe the signal process.

The defender is a Bayesian learner. Based on observations of the signal process $Y$ up to time $t$, the defender calculates \emph{suspicion level} $q_t$, the probability of the suspect being the attacker at time $t$:
\begin{equation}\label{bayesian}
q_t =   \PP \left(\theta=1  \big\vert  \FF_t^{Y} \right), 
\end{equation}
where $(\FF_t^{Y})_{t\geq0}$ is the filtration generated by the signal process $(Y_t)_{t \geq 0}$. We derive the SDE for $q_t$ that describes the change in the suspicion level given the attack intensity process $(\Delta_t)_{t \geq 0}$, using Theorem 8.1 of \cite{LiptserShiryaev}:
\begin{equation}
 \begin{split}\label{eqn-suspicion-level}
 dq_t&=\frac{1}{\sigma^2} \Big(\E[\theta \, \Delta_t 1_{\{  \theta=1\}}|\FF^Y_t] - \E[\theta|\FF^Y_t] \cdot \E[\Delta_t 1_{\{  \theta=1\}}|\FF^Y_t]\Big)\cdot\Big(dY_t - \E[\Delta_t 1_{\{  \theta=1\}}|\FF^Y_t] dt \Big)\\
&=\frac{q_t(1-q_t)\Delta_t}{\sigma^2} \Big(dY_t - q_t \Delta_t dt\Big).\\
\end{split}
 \end{equation}

\begin{figure}
	\centering
	\includegraphics[width=0.3\linewidth]{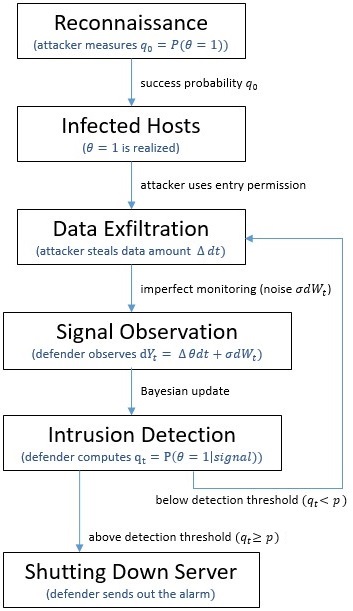}
	\caption{}
	\label{fig:gameflow}
\end{figure}

We assume that the game can be over by a random time $T$ that is independent of $\theta$ and $(W_t)_{t\geq 0}$ and has an exponential distribution
\begin{equation}\label{T_exp}
\PP(T>t)=e^{-rt},
\end{equation}
with a constant $r>0$. Note that $T$ is an exogenously given random variable that the players cannot control. The true value of $\theta$ is also publicly revealed if the game ends due to the random termination. One interpretation of the random termination time $T$ is that the suspect's identity $\theta$ can be revealed to the defender due to other independent factors.

The attacker's strategy $(\Delta_t)_{t \geq 0}$ is assumed to be a progressively measurable process with respect to the filtration $(\FF_t^Y)_{t \geq 0}$. The attacker obtains profits of $\Delta_t dt$ during infinitesimal time interval $dt$ until the game is over at $T\wedge \tau:=\min\{T,\tau\}$, where $\tau$ is the time the defender blocks the suspect. The attacker seeks the optimal attack intensity to maximize her expected profits:\footnote{If $\tau$ is a fixed random time, then the obvious optimal strategy in \eqref{AttackerProfit} is $\Delta_t=M$ for all $t\geq 0$. However, the defender bases his decision $\tau$ on observations of the signal process $Y$. Equations \eqref{eqn-suspicion-level} and \eqref{defenderCost} imply that $\tau$ will depend on the attack intensity $\Delta$, and thus the optimization in \eqref{AttackerProfit} is not obvious at all.}
\begin{equation}\label{AttackerProfit}
\max_{0 \leq (\Delta_t)_{t\geq0} \leq M} \E\left[\int_0^{T \wedge \tau} \Delta_t dt  \, \Big\vert \,  \theta=1 \right].
\end{equation}
As one can see from the expression above, the attacker's profit is larger if she steals a larger amount of data for a longer period of time.

If the suspect is the attacker, the defender incurs costs of $\Delta_t dt$ during infinitesimal time interval $dt$ until the defender blocks the suspect. If the suspect is innocent, the defender incurs zero costs during the game, but incurs a one-time {\it false alarm cost} $l > 0$ if the defender blocks the innocent user before the random termination. Even though we use the phrase `false alarm cost', $l$ does not necessarily represent penalty for wrong detection only. The term $l$ represents the defender's opportunity costs when he disables the innocent suspect. For instance, $l$ can include foregone profits during system downtimes, claims from customers because of inconvenience, and actual costs of inspecting the suspect's identity.

The defender's strategy is when to block the suspect, and his strategy is represented by a stopping time with respect to the filtration $(\FF_t^Y)_{t\geq0}$.  We denote $\mathcal{T}$ as the set of  all stopping times with respect to $(\FF_t^Y)_{t \geq 0}$. The defender's goal is to find the optimal stopping time to minimize expected costs:
\begin{equation}\label{defenderCost}
\min_{\tau \in \mathcal{T}}  \E\left[ 1_{\{\theta=1\}} \cdot \left(  \int_0^{T \wedge \tau}  \Delta_t   dt  \right)  + 1_{\{ \theta=0, \, \tau<T \}} \cdot  l  \right].
\end{equation}
The defender's cost is larger if a larger amount of data is stole or if the false alarm cost is larger, which seems reasonable to some extent.

In this paper, we restrict our attention to sequential equilibrium in Markov pure strategies that depend only on the suspicion level $q_t$ in \eqref{bayesian}. To be more specific, the attacker's equilibrium strategy is represented by a function $\alpha: [0,1] \to [0,M]$ of the suspicion level $q$, and the defender's equilibrium strategy is represented by a closed set $S\subset [0,1]$ which is the collection of suspicion levels at which the defender stops the game. For convenience, we denote $\tau_S$ as the first time the suspicion level process hits the closed set $S$,
\begin{equation}
\begin{split}\label{tau_S}
\tau_S = \inf \{ t\geq 0: q_t \in S \}.
\end{split}
\end{equation}
Clearly, $\tau_S$ is a stopping time with respect to $(\FF_t^Y)_{t\geq0}$.

Now, we introduce the definition of our Markov equilibrium.

\begin{definition}\label{def-equilibrium}
	Consider a process $(q_t)_{t\geq0}$, a closed set $S \subset [0,1]$, and a Lipschitz continuous function $\alpha: [0,1] \to [0,M]$. We say that the triplet $((q_t)_{t\geq0},S, \alpha)$ is a Markov equilibrium if the following conditions hold:
	
	(1) (Consistency) The process $(q_t)_{t \geq 0}$ satisfies Bayes' rule \eqref{bayesian}, given the initial value $q_0$ and the attack intensity $\Delta_t=\alpha(q_t)$.
		
	(2) (Attacker's optimality) The process $\big(\alpha(q_{t})\big)_{t \geq 0}$ is the solution to the attacker's profit maximization problem \eqref{AttackerProfit} for given $\tau=\tau_S$,
	\begin{equation}
	\begin{split}\label{attacker_max}
	\big(\alpha(q_{t})\big)_{t \geq 0} \in \argmax_{ 0\leq  (\Delta_t)_{t\geq 0} \leq M }  \E\left[   \int_0^{T\wedge \tau_S} \Delta_t   dt   \, \Big| \, \theta=1 \right].
	\end{split}
	\end{equation}

	(3) (Defender's optimality) The stopping time $\tau_S$ solves the defender's cost minimization problem \eqref{defenderCost} for given attack intensity $\Delta_t=\alpha(q_t)$,
	\begin{equation}
	\begin{split}\label{Dcost1}
	\tau_S \in \argmin_{\tau \in \mathcal{T}}  \E\left[ 1_{\{\theta=1\}} \cdot \left(  \int_0^{T \wedge \tau} \alpha(q_t)   dt  \right)  + 1_{\{ \theta=0, \, \tau<T \}} \cdot l  \right].
	\end{split}
	\end{equation}
\end{definition}

The first condition in the definition above implies that the suspicion level is calculated as if the initial value is $q_0$ and the attack intensity is $\alpha$, given the observation of the signal process. The second and third conditions are typical, meaning that every player's equilibrium strategy is the best response to everyone else's equilibrium strategy.

\section{Equilibrium Analysis}\label{equilibrium_analysis}

\subsection{Heuristic Derivation of Differential Equations When $S=[p,1]$}\label{subsec_derive}

We first present heuristic derivation of the differential equations that the players' value functions satisfy in Markov equilibrium $((q_t)_{t \geq 0}, [p,1], \alpha)$. Indeed, in \proref{S_unique_appendix}, we provide the result that the defender's equilibrium strategy should be the form of $S=[p,1]$. For simpler presentation, we slightly abuse our notation of stopping times as 
$$\tau_k=\inf{\{ t\geq0: q_t\in [k,1] \}} \quad \textrm{for} \quad k\in [0,1].$$

The attacker's expected profit can be written as
\begin{equation}
	\begin{split}\label{modified form attacker}
	\E\Big[\int_0^{T\wedge \tau_p}  \Delta_t dt \, \Big| \,  \theta=1\Big]&=\E\Big[\int_0^{\tau_p} 1_{\{T>t \}} \Delta_t dt \, \Big| \,  \theta=1\Big]=\E\Big[\int_0^{ \tau_p}  e^{-rt}\Delta_t dt \, \Big| \,  \theta=1\Big],
	\end{split}
	\end{equation}
where the second equality is from the independence of $T$ and the other random variables. Using the expression in \eqref{modified form attacker}, we define the value function $V$ as
\begin{equation}
\begin{split}\label{V_value}
V(q) :=&\max_{0 \leq (\Delta_t)_{t\geq0} \leq M}  \E\Big[\int_0^{\tau_p} e^{-rt}  \Delta_t dt   \, \Big| \,  \theta=1,\, q_0=q\Big].
\end{split}
\end{equation}
We derive the Hamilton-Jacobi-Bellman (HJB) equation for $V$ in \eqref{V_value}. The first condition in the definition of sequential equilibrium specifies how the suspicion level $q_t$ is calculated given the observation of the signal process up to time $t$. In our equilibrium concept, it is a common knowledge for both players that the attacker chooses attack intensities according to $\alpha$ in equilibrium. Therefore, when the attacker actually chooses $\Delta_t$ at time $t$, the infinitesimal change in the suspicion level is:
\begin{equation}\label{eqn-dq1}
dq_t = \frac{q_t(1-q_t) \alpha(q_t)}{\sigma^2} \left( \Delta_t dt + \sigma dW_t - q_t \alpha(q_t)dt  \right).
\end{equation}
Given equation \eqref{eqn-dq1}, the attacker's profit maximization problem \eqref{V_value} produces the following HJB equation:
\begin{equation}
\begin{split}\label{max}
&-r V(q) - V'(q) \tfrac{q^2(1-q) \alpha(q)^2}{\sigma^2} + \tfrac{1}{2} V''(q)  \tfrac{q^2(1-q)^2 \alpha(q)^2}{\sigma^2}   + \max_{\Delta\in [0,M]}   \big(V'(q) \tfrac{q(1-q) \alpha(q)}{\sigma^2} +1\big)\Delta  = 0,
\end{split}
\end{equation}
for $q\in (0,p)$. If $V'(q) \cdot\frac{q(1-q) \alpha(q)}{\sigma^2} +1 \geq  0$, then $\Delta = M$ maximizes the left hand side of \eqref{max}. Since the maximizer in \eqref{max} is supposed to be $\alpha(q)$ in equilibrium, we rewrite \eqref{max} for the case of $V'(q)  \geq -\tfrac{\sigma^2}{q(1-q)M}$,
\begin{equation}
\begin{split}\label{max_1}
\tfrac{V''(q)}{2}+ \tfrac{V'(q)}{q} - \tfrac{r\sigma^2 V(q)}{M^2 q^2 (1-q)^2} + \tfrac{\sigma^2}{M q^2 (1-q)^2}=0. \\
\end{split}
\end{equation}
For the case of $V'(q) < -\tfrac{\sigma^2}{q(1-q)M}$, we set $\alpha(q)=-\tfrac{\sigma^2}{q(1-q)} \cdot \frac{1}{V'(q)}$ to match the maximizer in \eqref{max} with $\alpha(q)$. Then any $\Delta\in [0,M]$ maximizes \eqref{max}, and we rewrite \eqref{max} as 
\begin{equation}
\begin{split}\label{max_2}
\tfrac{V''(q)}{2}- \tfrac{V'(q)}{1-q} - \tfrac{r}{\sigma^2} V'(q)^2 \,V(q)=0. \\
\end{split}
\end{equation}  
To find the boundary condition for the attacker's value function in \eqref{V_value}, we observe that $q_0=0$ implies $\tau_p=\infty$ because $0$ is an absorbing state (see SDE \eqref{eqn-suspicion-level}). Then the attacker chooses the highest intensity $M$ all the time and the corresponding value is $V(0)=\int_0^\infty e^{-rt} M dt = \tfrac{M}{r}$. 

Below are the differential equation and boundary conditions that the attacker's value function satisfies.
\begin{equation}
\begin{split}\label{hjb_V}
\begin{cases}
&\textrm{If  }  q\in (0,p) \textrm{  and  } V'(q) \geq -\tfrac{\sigma^2}{M(1-q)q}, \textrm{  then  }
\begin{cases} 
\frac{V''(q)}{2}+ \frac{V'(q)}{q} - \frac{r\sigma^2 V(q)}{M^2 q^2 (1-q)^2} + \frac{\sigma^2}{M q^2 (1-q)^2}=0,\\
\alpha(q)=M.
\end{cases}\\
&\textrm{If  }  q\in (0,p) \textrm{  and  } V'(q) < -\tfrac{\sigma^2}{M(1-q)q}, \textrm{  then  }
\begin{cases} 
\frac{V''(q)}{2}- \frac{V'(q)}{1-q} - \frac{r}{\sigma^2} V'(q)^2 \,V(q)=0,\\
\alpha(q) = -\frac{\sigma^2}{q(1-q)V'(q)}.
\end{cases}\\
&\textrm{If  }  q\in [p,1], \textrm{  then  } V(q)=0.\\
&V(0)= \frac{M}{r}.
\end{cases}
\end{split}
\end{equation}


Now we derive a differential equation from the defender's optimal stopping problem \eqref{Dcost1}. Using \eqref{bayesian}, the defender's expected cost at time $0$ can be written as 
\begin{equation}
\begin{split}
\E\Big[\int_0^\tau e^{-rt} \alpha(q_t) q_t \,dt + e^{-r \tau} l  (1-q_\tau) \Big].\\
\end{split}
\end{equation}
We define the value function $U$ as
\begin{equation}
\begin{split}\label{U_value}
U(q) :=\min_{\tau \in \mathcal{T}} \,  \E\Big[\int_0^\tau e^{-rt} \alpha(q_t) q_t \,dt + e^{-r \tau} l  (1-q_\tau)  \, \Big | \, q_0=q \Big].
\end{split}
\end{equation}

Since the defender's perception of the attack intensity is $\Delta_t=\alpha(q_t)$ in equilibrium, the corresponding SDE for the suspicion level is:
\begin{equation}\label{eqn-dq2}
dq_t = \frac{q_t(1-q_t) \alpha(q_t)}{\sigma^2} \left( \alpha(q_t) dt + \sigma dW_t - q_t \alpha(q_t)dt  \right).
\end{equation}
Given equation \eqref{eqn-dq2},
the optimal stopping problem \eqref{U_value} produces the following variational inequality and the characterization of the optimal stopping time $\tau_p$ in equilibrium:
\begin{align}
0&=\min \Big\{ -r U(q)  + \tfrac{1}{2} U''(q)  \cdot\tfrac{q^2(1-q)^2 \alpha(q)^2}{\sigma^2} + q \, \alpha(q)\, , \,\,  l(1-q) - U(q) \Big\}, \label{U_hjb_original}\\
\tau_p &= \inf\{ t \geq 0: U(q_t)=l(1-q_t) \}. \label{tau*}
\end{align}
If $q_0=0$, then $q_t \equiv 0$ all the time (see SDE \eqref{eqn-suspicion-level}). Therefore, in case $q_0=0$, the maximizer in \eqref{U_value} is $\tau=\infty$ and we obtain $U(0)=0$. With this boundary condition and the smooth-fit principle,
we rewrite the variational inequality \eqref{U_hjb_original} as  
\begin{equation}
\begin{split}\label{hjb_U}
\begin{cases}
&\textrm{If  }  q\in (0,p), \,\, \textrm{then  }
\left \{
\begin{tabular}{l}
$  -r U(q)  +   \tfrac{q^2(1-q)^2 \alpha(q)^2 U''(q)}{2\sigma^2} + q \, \alpha(q)=0$,   \\
$U(q) < l(1-q). $  
\end{tabular}
\right .  
\\
&\textrm{If  }  q\in [p,1], \,\, \textrm{then  }
\left \{
\begin{tabular}{l}
$  -r U(q)  +   \tfrac{q^2(1-q)^2 \alpha(q)^2 U''(q)}{2\sigma^2} + q \, \alpha(q)\geq0$,   \\
$U(q) =l(1-q). $  
\end{tabular}
\right .  
\\
&U(0)=0, \quad  \lim_{q\uparrow p}U(q)=l(1-p), \quad \lim_{q\uparrow p}U'(q)=-l. 
\end{cases}
\end{split}
\end{equation}
In summary, \eqref{hjb_V} and \eqref{hjb_U} constitute the system of differential equations for the Markov equilibrium. To describe the explicit solution of the system, we first define the functions $\varphi, y$ and constants $a, b, c, q^*$ as 
\begin{equation}
\begin{split} \label{const}
\varphi(x)&:=\tfrac{2}{\sqrt{\pi}}\int_0^x e^{-t^2} dt, \\
 a&:=\tfrac{1}{2}\Big( \sqrt{1+\tfrac{8r\sigma^2}{M^2}} - 1  \Big), \\ 
 b&:=\tfrac{(1-a)M}{2\sigma\sqrt{r}},\\
c&:=\tfrac{2\sigma\sqrt{r}e^{-b^2} }{M  \sqrt{\pi}} +  \varphi(b),\\
  q^*&:=\tfrac{p(c-\varphi(b) )}{c-p \, \varphi(b) },\\
  y(x)&:= \varphi^{-1}\big(\tfrac{c (p - x)}{p (1 - x)} \big).
\end{split}
\end{equation}
Note that the constant $q^*$ and the function $y$ depend on the constant $p$ that will be determined later. 

\begin{lemma}
	\label{tech_lemma}
	Let the functions $\varphi, y$ and constants $a, b, c, q^*$ be as in \eqref{const}. If $\tfrac{r \sigma^2}{M^2} < 1$, then the followings hold.
	
	(1) $a\in (0,1)$, $b>0$, $c\geq 1$, $q^*\in (0,p)$ and $y(q^*)=b$.
	
	
	(2) $0<c\sqrt{\pi}  \, \big( \tfrac{(1-p)x}{p (1 - x)}\big) e^{y(x)^2} -  \tfrac{2\sigma \sqrt{r}}{M} $ for $x\in (q^*,p]$.
\end{lemma}
\begin{proof}
	Elementary calculations produce (1). To obtain (2), we first observe that $z\mapsto \varphi(z)+\tfrac{2 \sigma \sqrt{r} \, e^{-z^2}}{M\sqrt{\pi}}$ is a strictly increasing function on $z\in [0,b]$. Then, for $x\in (q^*,p]$, we observe that $y(x)< b$ and
	$$\varphi(b)+\tfrac{2 \sigma \sqrt{r} \, e^{-b^2}}{M\sqrt{\pi}} > \varphi(y(x))+\tfrac{2 \sigma \sqrt{r} \, e^{-y(x)^2}}{M\sqrt{\pi}} = \Big(\tfrac{2 \sigma \sqrt{r} \, e^{-y(x)^2}}{M\sqrt{\pi}} -\tfrac{c(1-p)x}{p (1 - x)} \Big) + c.$$
	The above inequality, together with the expression of $c$, produces (2).
\end{proof}

The unique explicit solution of the system \eqref{hjb_V} and \eqref{hjb_U} is provided in the next proposition.

\begin{proposition}\label{solutions} 
	The unique solution of the system \eqref{hjb_V} and \eqref{hjb_U} that satisfies $V,U \in C^2([0,p))$ has the following expression:
	
	\noindent(1) If $\tfrac{r \sigma^2}{M^2} \geq  1$, then 
	\begin{align}
	p&=\tfrac{(1+a)rl}{(1+a)rl + a M }\\
	 \alpha(q)&=M \quad \textrm{for} \quad q\in [0,1]\\
	V(q)&=\begin{cases}
	\tfrac{M}{r}\big(1- \big(\tfrac{1-p}{p}\big)^a \big(\tfrac{q}{1-q}\big)^a \big) & \textrm{for  }\,\, q\in [0,p)\\
	0 & \textrm{for  }\,\, q \in [p,1]
	\end{cases}  \label{V_exp1}\\
	U(q)&=\begin{cases}
	q \Big( \tfrac{M}{r} - \big( \tfrac{M}{r} - \tfrac{(1-p)l}{p} \big)\big(\tfrac{1-p}{p}\big)^a \big(\tfrac{q}{1-q}\big)^a    \Big)  &\textrm{for  }\,\, q\in [0,p)\\
	(1-q)l &\textrm{for  }\,\, q \in [p,1]
	\end{cases} \label{U_exp1}
	\end{align}
	\noindent	(2) If $\tfrac{r \sigma^2}{M^2} <  1$, then
	\begin{align}
	p&=\tfrac{c \, l \sqrt{\pi r}}{c \, l \sqrt{\pi r} + \sigma  }\\
	\alpha(q)&=\begin{cases}M &\textrm{for  }\,\, q \in [0,q^*], \\
\frac{2p(1-q)\sigma \sqrt{r} }{c \sqrt{\pi}\,(1-p)q}\,e^{-y(q)^2} &\textrm{for  }\,\,  q\in (q^*,p),\\
\tfrac{2\sigma \sqrt{r}}{c\sqrt{\pi}} &\textrm{for  }\,\, q \in [p,1].
\end{cases}\\
	V(q)&=\begin{cases}\frac{M}{r} - \frac{\sigma^2}{aM} \big(\frac{1-q^*}{q^*}\big)^a  \big(\frac{q}{1-q}\big)^a &\textrm{for  }\,\, q \in [0,q^*] \\
	\frac{\sigma }{\sqrt{r}}\,y(q) &\textrm{for  }\,\,  q\in (q^*,p)\\
	0 &\textrm{for  }\,\, q \in [p,1]
	\end{cases}  \label{V_exp2}\\
	U(q)&=\begin{cases}
	q\Big(\tfrac{M}{r} - \big( \tfrac{\sigma^2}{aM} - \frac{c\sqrt{\pi r}  (1-p) l \sigma}{(1+a)M p } \big)   \big(\frac{1-q^*}{q^*}\big)^a  \big(\frac{q}{1-q}\big)^a  \Big)  &\textrm{for  }\,\, q \in [0,q^*] \\
	\frac{\sigma }{\sqrt{r}}\,q\,y(q)  +  l(1-q)   \big( e^{-y(q)^2} -  \tfrac{c \sqrt{\pi}(1-p)q\,y(q)}{p(1-q)}    \big)   &\textrm{for  }\,\,  q\in (q^*,p)\\
	(1-q)l &\textrm{for  }\,\, q \in [p,1]
	\end{cases}   \label{U_exp2}
	\end{align}	
\end{proposition}

\begin{proof}
	\lemref{tech_lemma} and explicit computations produce the proposition.
\end{proof}

\subsection{Characterization of the Unique Equilibrium}

We further restrict our attention to Markov equilibria in which the value functions of the attacker and the defender are smooth enough (twice differentiable), and show that there exists a unique Markov equilibrium that induces the smooth enough value functions.

\begin{theorem}\label{main_thm}
There exists a unique\footnote{It turns out that in the Markov equilibrium $((q_t)_{t\geq0},S,\alpha)$, $\alpha \vert_{(p,1]}$ has little to no impact on the stopping threshold $p$. That is, for another Lipschitz continuous function $\hat \alpha$ such that $\hat \alpha \vert_{[0,p]} = \alpha \vert_{[0,p]}$, the defender optimally chooses $p$ as long as $\hat\alpha \vert_{(p,1]}$ is not too small. This means that $((q_t)_{t\geq0},S, \hat \alpha)$ is another Markov equilibrium. However, the pair of the optimal stopping threshold $p$ and the optimal attack intensity $\alpha$ over the interval $[0,p]$ is uniquely determined by the exogenous parameters. In this sense we say that our Markov equilibrium is unique. For simplicity, we set $\alpha(q)=\alpha(p)$ for $q\in (p,1]$.
} Markov equilibrium $((q_t)_{t\geq0},S, \alpha)$. And, there exists a unique $p \in [0,1]$ such that $S =  [ p,1]$.
The equilibrium stopping threshold $p$ and the equilibrium attack intensity $\alpha$ have the following form:
\begin{equation}\label{p alpha form}
\begin{split}
\textrm{In case}\quad  \tfrac{r \sigma^2}{M^2}\geq 1:&\quad p=\tfrac{(1+a)rl}{(1+a)rl + a M },\quad  \alpha(q)=M \quad \textrm{for  }q\in [0,1]. \\
\textrm{In case}\quad \tfrac{r \sigma^2}{M^2}< 1:& \quad p=\tfrac{c \, l \sqrt{\pi r}}{c \, l \sqrt{\pi r} + \sigma  },\quad \alpha(q)=\begin{cases}M &\textrm{for  }\,\, q \in [0,q^*], \\
\frac{2p(1-q)\sigma \sqrt{r} }{c \sqrt{\pi}\,(1-p)q}\,e^{-y(q)^2} &\textrm{for  }\,\,  q\in (q^*,p),\\
\tfrac{2\sigma \sqrt{r}}{c\sqrt{\pi}} &\textrm{for  }\,\, q \in [p,1].
\end{cases}
\end{split}
\end{equation}
\end{theorem}
\begin{proof}
We postpone the proof to Section 4. The result is the direct consequence of \proref{main_appendix} and \proref{S_unique_appendix}.
\end{proof}
\noindent
This theorem shows that in equilibrium, the defender blocks the suspect once the suspicion level exceeds a certain threshold $p$, allowing us to use the {\it stopping threshold} instead of the set of suspicion levels at which the defender blocks the suspect.

The expression of the equilibrium in \eqref{p alpha form} implies that if  $\tfrac{r \sigma^2}{M^2} \geq 1$, it is the dominant strategy for the attacker to choose the highest attack intensity $M$ all the time. To understand this, it is helpful to imagine three extreme cases in which  $\tfrac{r \sigma^2}{M^2}$ is very large. A high probability $r$ of random termination implies that the game is more likely to end due to random termination, which in turn implies that the attacker has a weaker incentive to slow down the defender's learning. A large noise $\sigma$ enables the attacker to hide behind the noise. If the upper bound $M$ of attack intensity is low, then the stopping threshold would be high because the aggregate running costs are low compared to the false alarm costs. In these cases, the attacker has an incentive to set the highest attack intensity.
If $\tfrac{r \sigma^2}{M^2}< 1$, the attacker chooses the highest attack intensity $M$ when the suspicion level $q_t$ is sufficiently low ($q_t \leq q^*$). As the suspicion level $q_t$ increases above $q^*$, the attacker gradually decreases the attack intensity to lower the rate at which the suspicion level is updated.

\begin{remark}\label{remark1}
Let us discuss the stopping threshold. Intuitively, if $p$ is the equilibrium stopping threshold, the defender is indifferent between stopping the game and continuing the game when the suspicion level is $p$. The defender incurs the expected cost of $(1-p)l$ if she stops the game at $p$. If the defender waits until the suspicion level becomes $p+dp$, the expected running cost increases and the expected false alarm cost decreases. The defender incurs aggregate running costs that she would not have incurred if she had stopped the game immediately. However, since the defender stops the game at a higher threshold, the expected false alarm cost decreases. Up to the first order of $dp$, it should be true that:
\begin{equation}\label{eqn-mcmb}
	\mathbb{E} \left[ \int_0^{\tau_{p+dp}\wedge T}  \alpha(q_t) dt \cdot 1_{\{ \theta = 1 \}}  + l \cdot 1_{\{\tau_{p+dp} < T, \theta = 0\}} \, \Big| \, q_0=p \right] = \mathbb{E}[l \cdot 1_{\{ \theta = 0  \}} \,|\, q_0=p].
\end{equation}
Rearranging this equation, up to the first order of $dp$, we obtain:
\begin{equation}\label{eqn-mbmc2}
		\mathbb{E} \left[ \int_0^{\tau_{p+dp}} e^{-rt}  \alpha(q_t) dt \cdot 1_{\{ \theta = 1 \}}  \, \Big| \, q_0=p   \right] =  \mathbb{E}\left[ \left( \int_0^{\tau_{p+dp}} r \, e^{-rt}       dt \right) l \cdot  1_{\{ \theta = 0  \}}  \, \Big| \, q_0=p  \right].
\end{equation}
This equation makes sense because it basically implies that marginal benefit equals marginal cost. The left hand side of the equation is the marginal increase in the expected cost and the right hand side is the marginal reduction in the false alarm cost.
\end{remark}
\begin{remark}
Based on the arguments in Remark \ref{remark1}, we can infer a possible impact of other types of running costs on the equilibrium. For instance, let $Z$ be a constant that represents a running cost such as monitoring cost. To be specific, we add the term $\E\left[\int_0^{T \wedge \tau} Z dt \right]$ to the defender's cost.

Similar to equation \eqref{eqn-mbmc2}, up to the first order of $dp$, we obtain:
\begin{equation}\label{eqn-mbmc3}
	\mathbb{E} \left[ \int_0^{\tau_{p+dp}} e^{-rt} \big(  \alpha(q_t) + Z \big)  dt \cdot 1_{\{ \theta = 1 \}}   \, \Big| \, q_0=p  \right] =  \mathbb{E}\left[ \int_0^{\tau_{p+dp}} e^{-rt}  ( r l - Z ) dt \cdot  1_{\{ \theta = 0  \}}   \, \Big| \, q_0=p \right].
\end{equation}
As we can see from the equation above, other types of running costs decrease the stopping threshold. Intuitively, if the defender incurs a larger amount of running costs (due to the monitoring cost $Z$), he is more willing to stop the game earlier to save costs.\footnote{Note that $Z$ does not have to be interpreted as costs only. The term $Z$ can be considered as income flow if it is negative, in which case the stopping threshold increases. If the defender earns positive profits, he is willing to take more risk of cyberattacks and tries to stop the game later.}

Equation \eqref{eqn-mbmc3} has another important implication. As one can infer from the previous paragraph, other types of running costs or income flows do change quantitative properties of equilibrium, but do not alter qualitative properties of equilibrium. To be more specific, in the model with the additional cost term $Z$, one can check that \thmref{main_thm} still holds with a different expression of $p$ that depends on $Z$.
\end{remark}

The explicit expression in \eqref{p alpha form} allows us to describe how the equilibrium threshold and attack intensity change as the exogenous parameter changes. 

\begin{proposition}\label{coro-alpha-comp-statics}
(1) The equilibrium stopping threshold $p$ increases in $l$ and $r$ and decreases in $\sigma$ and $M$.

(2) The equilibrium attack intensity  $\alpha$ increases in $l,r,\sigma$ and $M$.

\end{proposition}
\begin{proof}
(1) We prove that $p$ is a decreasing function of $M$. Other cases can be proved similarly. In case $\tfrac{r\sigma^2}{M^2}\geq 1$, we have $p=\tfrac{(1+a)rl}{(1+a)r l + aM}$. We substitute $a$ in \eqref{const} into the expression of $p$ and compute the derivative,
$$
\tfrac{\partial p}{\partial M}=-\tfrac{8l r^2\sigma^2 (\sqrt{M^2+8r\sigma^2}-2M)}{\sqrt{M^2+8r\sigma^2}\big( -M^2+l r\sqrt{M^2+8r\sigma^2}+M(l r+\sqrt{M^2+8r\sigma^2}) \big)^2}<0,
$$
where the inequality is due to $\tfrac{r\sigma^2}{M^2}\geq 1$.

In case $\tfrac{r\sigma^2}{M^2}< 1$, we have $p=\tfrac{cl \sqrt{\pi r}}{cl \sqrt{\pi r}+\sigma}$. We substitute $c$ in \eqref{const} into the expression of $p$ and compute the derivative,
$$
\tfrac{\partial p}{\partial M} =-\tfrac{e^{-b^2}l \big((M^2+2r\sigma^2) \sqrt{M^2+8r\sigma^2}-M^3-6Mr\sigma^2  \big)}{M^2\sqrt{M^2+8r\sigma^2}(\sigma+l \sqrt{\pi r}c)^2}<0,
$$
where the inequality is due to $(M^2+2r\sigma^2)^2(M^2+8r\sigma^2)-(M^3+6Mr\sigma^2)^2=32r^3\sigma^6>0 $. 

\smallskip

(2) It is enough to check the monotonicity of $\alpha(q)$ in the parameters, for $\tfrac{r\sigma^2}{M^2}< 1$ and $q\in (q^*,p]$. We substitute $p$ in \eqref{p alpha form} to the expression of $\alpha$ in \eqref{p alpha form} and obtain
\begin{align}\label{alpha_proof1}
\alpha(q)=\tfrac{2(1-q)}{q}\, l \, r \, \exp{\big(-\varphi^{-1}(c- \tfrac{q \sigma}{\sqrt{\pi r}l(1-q)})^2\big)}\quad \textrm{for}\quad q\in (q^*,p].
\end{align}

(i) ($\alpha$ increases in $l$): As we know that $p$ increases in $l$ by part (1), the expression of $\alpha$ in \eqref{p alpha form} for $q\in (q^*,p]$ implies that it is enough to check that $\alpha$ increases in $p$. We first check that 
\begin{align}\label{ineq_1}
e^{-y(q)^2}-c\sqrt{\pi}(\tfrac{(1-p)q}{p(1-q)}) y(q) \geq 0 \quad \textrm{for  } q\in (q^*,p],
\end{align}
where $y(q)$ is defined in \eqref{const}. Indeed, we observe that
\begin{displaymath}
\begin{split}
&(e^{-y(q)^2}-c\sqrt{\pi}(\tfrac{(1-p)q}{p(1-q)}) y(q)\big)\big|_{q=q^*}=ae^{-b^2}>0,\\
&\tfrac{d}{dq}\big(e^{-y(q)^2}-c\sqrt{\pi}(\tfrac{(1-p)q}{p(1-q)}) y(q)\big)=\tfrac{c^2\pi q (1-p)^2 e^{y(q)^2}}{2p^2(1-q)^3}>0,
\end{split}
\end{displaymath}
and conclude the inequality \eqref{ineq_1}. Then, we observe that $\alpha$ increases in $p$ due to \eqref{ineq_1} and the following expression: 
\begin{align}
\tfrac{\partial}{\partial p}\alpha(q)=\tfrac{2\sqrt{r}\sigma(1-q)}{c\sqrt{\pi}(1-p)^2 q}\big(  e^{-y(q)^2}-c\sqrt{\pi}\big(\tfrac{(1-p)q}{p(1-q)}\big) y(q)\big).
\end{align}

(ii) ($\alpha$ increases in $\sigma$): Due to the expression of $\alpha$ in \eqref{alpha_proof1}, it is enough to show that \\$c- \tfrac{q \sigma}{\sqrt{\pi r}l(1-q)}$ decreases in $\sigma$, because the map $x\mapsto e^{-\varphi^{-1}(x)^2}$ decreases in $x$ for $x>0$. We observe that for $q\in (q^*,p]$,
\begin{displaymath}
\begin{split}
\tfrac{\partial}{\partial \sigma} \big(  c- \tfrac{q \sigma}{\sqrt{\pi r}l(1-q)} \big) &= \tfrac{\partial c}{\partial \sigma} - \tfrac{q}{\sqrt{\pi r}l(1-q)} <   \tfrac{\partial c}{\partial \sigma} - \tfrac{q^*}{\sqrt{\pi r}l(1-q^*)}=\tfrac{e^{-b^2} (M\sqrt{M^2+8 r \sigma^2} - M^2 - 6 r \sigma^2) }{\sigma^2\sqrt{\pi r (M^2+8 r \sigma^2)}} <0.
\end{split}
\end{displaymath}
The last inequality above is from $M^2(M^2+8 r \sigma^2)-(M^2 + 6 r \sigma^2)^2=-4r \sigma^2(M^2+9 r \sigma^2)<0$.

(iii) ($\alpha$ increases in $r$): We first observe that
\begin{align}\label{c_r}
\tfrac{\partial c}{\partial r} = \tfrac{((2A^2+1)\sqrt{1+8A^2}-(6A^2+1))e^{-b^2}}{2\sqrt{\pi}r A \sqrt{1+8A^2}}>0, \quad \textrm{where}\quad A=\tfrac{\sigma \sqrt{r}}{M}.
\end{align}
The inequality above is from $((2A^2+1)\sqrt{1+8A^2})^2- (6A^2+1)^2=32A^6$. We also observe
\begin{equation}
\begin{split}\label{c_r2}
&\tfrac{\partial}{\partial q}\Big(e^{-y(q)^2}- \sqrt{\pi} y(q)\big(\tfrac{c(1-p)q}{2p(1-q)}+r \tfrac{\partial c}{\partial r} \big)\Big)\\
&= \tfrac{c(1-p)\sqrt{\pi} \big(2p(1-q)y(q)+e^{y(q)^2} \sqrt{\pi} (c(1-p)q+2pr(1-q)  \tfrac{\partial c}{\partial r}) \big)}{4p^2(1-q)^3}>0, \quad \textrm{for  }q\in (q^*,p],
\end{split}
\end{equation}
where we use \eqref{c_r} for the inequality. Using the expression of $\alpha$ in \eqref{alpha_proof1} and the inequalities \eqref{c_r} and \eqref{c_r2}, we observe that for $q\in (q^*,p]$,
\begin{equation}
\begin{split}
\tfrac{\partial}{\partial r} \alpha(q)&=\tfrac{2(1-q)l}{q} \Big(e^{-y(q)^2}- \sqrt{\pi} y(q)\big(\tfrac{c(1-p)q}{2p(1-q)}+r \tfrac{\partial c}{\partial r} \big)\Big)\\
&>\tfrac{2(1-q)l}{q} \Big(e^{-y(q^*)^2}- \sqrt{\pi} y(q^*)\big(\tfrac{c(1-p)q^*}{2p(1-q^*)}+r \tfrac{\partial c}{\partial r} \big)\Big)\\
&=\tfrac{2(1-q)l e^{-b^2}\big( 16A^4+15A^2+2 - (5A^2+2)\sqrt{1+8A^2})  \big)}{4q A^2 \sqrt{1+8A^2}}, \quad \textrm{where}\quad A=\tfrac{\sigma \sqrt{r}}{M}.
\end{split}
\end{equation}
The above expression and the inequality $(16A^4+15A^2+2)^2-((5A^2+2)\sqrt{1+8A^2})^2=8(32A^8+35A^6+13A^4+A^2)>0$
produce $\frac{\partial}{\partial r} \alpha(q)>0$.

(iv) ($\alpha$ increases in $M$): Due to the expression of $\alpha$ in \eqref{alpha_proof1}, it is enough to show that $c$ decreases in $M$, because the map $x\mapsto e^{-\varphi^{-1}(x)^2}$ decreases in $x$ for $x>0$. Indeed, 
\begin{align}
\tfrac{\partial c}{\partial M} = \tfrac{e^{-b^2}(M^3+6Mr \sigma^2 - (M^2+2r\sigma^2)\sqrt{M^2+8r\sigma^2})}{\sqrt{\pi r} M^2 \sigma \sqrt{M^2 +8r \sigma^2}}<0,     \nonumber
\end{align}
where the last inequality is by $(M^3+6Mr \sigma^2)^2 - ((M^2+2r\sigma^2)\sqrt{M^2+8r\sigma^2}))^2=-32r^3\sigma^6<0$. 
\end{proof}

\section{Proof of \thmref{main_thm}}

This section is devoted to the proof of \thmref{main_thm}. In \proref{solutions}, we provide the unique solution of the system of the differential equations \eqref{hjb_V} and \eqref{hjb_U}. In \proref{main_appendix}, we verify that the unique solution indeed constitutes a Markov equilibrium. In \proref{S_unique_appendix}, we show that if $((q_t)_{t\geq0},S,\alpha)$ is a Markov equilibrium, then the set $S$ should be of the form $[p,1]$ for a constant $p>0$. All in all, \proref{main_appendix} and \proref{S_unique_appendix} complete the proof of \thmref{main_thm}.

\begin{proposition}\label{main_appendix}
	$((q_t)_{t\geq 0},[p,1],\alpha)$ defined in \eqref{p alpha form} is a Markov equilibrium in \defref{def-equilibrium}, and $V$ and $U$ in \proref{solutions} are the value functions of the attacker and defender.
\end{proposition}
\begin{proof}

	{\bf Checking (1) in \defref{def-equilibrium}}
	
	The Lipschitz continuity of $\alpha$ ensures that the SDE \eqref{eqn-suspicion-level} has a unique solution when $\Delta_t=\alpha(q_t)$. Theorem 8.1 of \cite{LiptserShiryaev} implies that the solution of the SDE \eqref{eqn-suspicion-level} satisfies \eqref{bayesian}.
	
	\medskip
	
	{\bf Checking (2) in \defref{def-equilibrium}}
	
	We prove the optimality of $\tau_p$ in \eqref{Dcost1}. For any $\tau\in \mathcal{T}$ and $U$ in \proref{solutions}, Ito's formula produces
	\begin{equation}
	\begin{split}\label{verify_Y1}
	&e^{-r (t \wedge \tau)} U(q_{t \wedge \tau}) + \int_0^{t\wedge \tau} e^{-r s} \alpha(q_s) q_s \, ds\\
	&= U(q_0) + \int_0^{t\wedge \tau} e^{-r s} \Big(-r U(q_s)  +   \tfrac{q_s^2(1-q_s)^2 \alpha(q_s)^2 }{2\sigma^2}  U''(q_s)+ q_s \, \alpha(q_s)        \Big) ds + \int_0^{t\wedge \tau} e^{-r s} U'(q_s) dq_s\\
	& \geq  U(q_0) +\int_0^{t\wedge \tau} e^{-r s} U'(q_s) dq_s,
	\end{split}
	\end{equation}
	where the inequality is due to the fact that $U$ satisfies \eqref{hjb_U}. If we consider the stopping time $\tau_p$, the inequality becomes an equality:
	\begin{equation}
	\begin{split}\label{verify_Y2}
	&e^{-r (t \wedge \tau_p)} U(q_{t \wedge \tau_p}) + \int_0^{t\wedge \tau_p} e^{-r s} \alpha(q_s) q_s \, ds= U(q_0) +\int_0^{t\wedge \tau_p} e^{-r s} U'(q_s) dq_s.
	\end{split}
	\end{equation}
	We apply the Fubini's theorem and the iterated conditioning to obtain
	\begin{equation}
	\begin{split}\label{verify_Y3}
	\E\Big[\int_0^{t\wedge \tau} e^{-r s} U'(q_s) dq_s \Big]&=\int_0^{\infty} \E\Big[  1_{\{ 0\leq s \leq t \wedge \tau  \}} e^{-r s}U'(q_s) \tfrac{q_s(1-q_s)\alpha(q_s)^2 }{\sigma^2}\E\big[1_{\{ \theta=1 \}} - q_s  \big| \FF_s^Y \big]     \Big] ds \\
	&\qquad + \E\Big[ \int_0^{t\wedge \tau} e^{-r s}U'(q_s) \tfrac{q_s(1-q_s)\alpha(q_s) }{\sigma} dW_s  \Big]\\
	&=0,
	\end{split}
	\end{equation}
	where the last equality holds because  $\E[1_{\{ \theta=1 \}} - q_s  \big| \FF_s^Y]=0$ (see \eqref{bayesian}) and the stochastic integral part is a square-integrable martingale with respect to the filtration $(\FF_t^Y)_{t \geq 0}$.
	We combine \eqref{verify_Y1}-\eqref{verify_Y3} and obtain
	\begin{equation}
	\begin{split}\label{verify_Y4}
	&\E\Big[e^{-r (t \wedge \tau)} U(q_{t \wedge \tau}) + \int_0^{t\wedge \tau} e^{-r s} \alpha(q_s) q_s \, ds\Big]\geq U(q_0)=\E\Big[e^{-r (t \wedge \tau_p)} U(q_{t \wedge \tau_p}) + \int_0^{t\wedge \tau_p} e^{-r s} \alpha(q_s) q_s \, ds\Big].
	\end{split}
	\end{equation}
	Since $U$ and $\alpha$ are bounded, as $t\to \infty$, the dominated convergence theorem produces
	\begin{align}
	\E\Big[e^{-r \tau} U(q_{\tau}) + \int_0^{\tau} e^{-r s} \alpha(q_s) q_s \, ds\Big]\geq \E\Big[e^{-r \tau_p} U(q_{ \tau_p}) + \int_0^{ \tau_p} e^{-r s} \alpha(q_s) q_s \, ds\Big]. \nonumber
	\end{align}
	The above inequality, together with \eqref{hjb_U}, implies that
	\begin{equation}
	\begin{split}\label{verify_Y6}
	&\E\Big[e^{-r \tau} l (1-q_\tau) + \int_0^{\tau} e^{-r s} \alpha(q_s) q_s \, ds\Big] \geq \E\Big[e^{-r \tau_p} l(1-q_{\tau_p}) + \int_0^{ \tau_p} e^{-r s} \alpha(q_s) q_s \, ds\Big].\\
	\end{split}
	\end{equation}
	It remains to derive \eqref{Dcost1} from \eqref{verify_Y6}. Since the process $(q_t)_{t\geq0}$ is uniformly bounded, we apply the optional sampling theorem\footnote{See, for example, \cite{KaratzasShreve} Theorem 3.22 in Chapter 1.} to the martingale $q_t= \E[1_{\{\theta=1\}}|\FF_t^Y]$ and obtain
	\begin{align}\label{verify_Y7}
	q_\tau  = \E[1_{\{\theta=1\}} | \FF_\tau^Y] \quad \textrm{for all  }\tau\in \mathcal{T}.
	\end{align}
	Using \eqref{verify_Y7} and the independence of $T$, we obtain the following equalities:
	\begin{equation}
	\begin{split}\label{verify_Y9}
	\E\Big[e^{-r \tau} l(1-q_\tau) + \int_0^{\tau} e^{-r s} \alpha(q_s) q_s \, ds\Big]&=\E\Big[   l \cdot e^{-r \tau}  1_{\{\theta=0\}} \Big]+  \int_0^{\infty}\E \Big[ 1_{\{s<\tau \}}e^{-r s} \alpha(q_s)  1_{\{\theta=1\}} \Big]\, ds\\
	&=\E\Big[   l \cdot  1_{\{\theta=0\}}  1_{\{ T >\tau \}}  \Big]+\int_0^{\infty}\E \Big[ 1_{\{s<\tau \}} 1_{\{ s<T \}}  \alpha(q_s)  1_{\{\theta=1\}} \Big]\, ds\\
	&= \E\Big[ \Big(\int_0^{T\wedge \tau} \alpha(q_s)  ds \Big) \cdot 1_{\{\theta=1\}} +  l \cdot  1_{\{\theta=0, \, \tau<T\}} \Big],
	\end{split}
	\end{equation}
	where we apply the Fubini's theorem for the first and third equality, and use the iterated conditioning for the first and second equality.
	Since \eqref{verify_Y6} and \eqref{verify_Y9} hold for any $\tau\in \mathcal{T}$, we conclude that $\tau_p$ is optimal in \eqref{Dcost1}.
	
	\medskip

	{\bf Checking (3) in \defref{def-equilibrium}}\\
	In this part of the proof, we use notation $q_t^{(\Delta)}$ and $\tau^{(\Delta)}_p$ instead of $q_t$ and $\tau_p$, to emphasize their dependence on the attacker's (possibly off-equilibrium) strategy $\Delta$. 
	To be specific, for attack intensity process $(\Delta_t)_{t\geq0}$, let the process $(q_t^{(\Delta)})_{t\geq0}$ be the solution of SDE \eqref{eqn-suspicion-level} and $\tau^{(\Delta)}_p=\inf\{t\geq 0:  q^{(\Delta)}_t \geq p\}$.\\
	To verify that the function $V$ in \proref{solutions} is indeed the optimal value of the attacker, we apply Ito's formula, conditioned on $\theta=1$: For $q_0\in [0,p]$,
	\begin{equation}
	\begin{split}\label{verify}
	& e^{-r(t\wedge \tau^{(\Delta)}_p)} V(q_{t\wedge \tau^{(\Delta)}_p}^{(\Delta)})+ \int_0^{t\wedge \tau^{(\Delta)}_p} e^{-rs} \Delta_s ds\\
	&= V(q_0)+\int_0^{t\wedge \tau^{(\Delta)}_p} e^{-rs} \Big( -r V(q) - \tfrac{q^2(1-q) \alpha(q)^2 V'(q) }{\sigma^2} + \tfrac{q^2(1-q)^2 \alpha(q)^2 V''(q)}{2\sigma^2} \\
	& \qquad  +   \big(\tfrac{q(1-q) \alpha(q)V'(q)}{\sigma^2} +1\big)\Delta_s \Big)\Big|_{q=q_s^{(\Delta)}} ds 
	+ \int_0^{t\wedge \tau^\Delta_p}  e^{-rs} \, \tfrac{q(1-q) \alpha(q)V'(q)}{\sigma}\Big|_{q=q_s^{(\Delta)}}  \, dW_s\\
	&\leq V(q_0) + \int_0^{t\wedge \tau^\Delta_p}  e^{-rs} \, \tfrac{q(1-q) \alpha(q)V'(q)}{\sigma}\Big|_{q=q_s^{(\Delta)}}   dW_s,
	\end{split}
	\end{equation}
	where the inequality above is due to \eqref{max}. Indeed, $V$ in \proref{solutions} satisfies \eqref{hjb_V}, and \eqref{hjb_V} implies \eqref{max}. Since $ \tfrac{q(1-q) \alpha(q)V'(q)}{\sigma}$ is bounded on $q\in[0,p)$, the stochastic integral term in \eqref{verify} is a square-integrable martingale (with respect to $(\FF_t^Y)_{t \geq 0}$) and has mean zero. Since the maximum is achieved at $\Delta = \alpha(q)$ in \eqref{max}, the inequality \eqref{verify} implies
	\begin{equation}
	\begin{split}\label{verify2}
	&\E \Big[ e^{-r(t\wedge \tau^{(\Delta)}_p)} V(q_{t\wedge \tau^{(\Delta)}_p}^{(\Delta)})+ \int_0^{t\wedge \tau^{(\Delta)}_p} e^{-rs} \Delta_s ds \, \Big| \, \theta=1\Big] \\
	&\leq V(q_0)=\E \Big[ e^{-r(t\wedge \tau_p^{(\alpha)})} V(q_{t\wedge \tau_p^{(\alpha)}}^{(\alpha)})+ \int_0^{t\wedge \tau_p^{(\alpha)}} e^{-rs} \alpha(q_s^{(\alpha)}) ds  \, \Big| \, \theta=1 \Big], 
	\end{split}
	\end{equation}
	where we denote $q_t^{(\alpha)}$ and $\tau_p^{(\alpha)}$ as the suspicion level process and the stopping time with $\Delta_t=\alpha(q_t^{(\alpha)})$. In \eqref{verify2}, we let $t\to \infty$, and the boundedness of $V$ produces  
	\begin{equation}
	\begin{split}\label{verify3}
	&\E\Big[ e^{-r\tau^{(\Delta)}_p} V(q_{\tau^{(\Delta)}_p}) + \int_0^{\tau^{(\Delta)}_p} e^{-rs} \Delta_s ds  \, \Big| \, \theta=1 \Big]\\
	& \leq V(q_0)=\E \Big[e^{-r\tau_p^{(\alpha)}} V(q_{\tau_p}^{(\alpha)}) + \int_0^{\tau_p^{(\alpha)}} e^{-rs} \alpha(q_s^{(\alpha)}) ds  \, \Big| \, \theta=1 \Big]. 
	\end{split}
	\end{equation}
	$V(p)=0$ implies that $e^{-r\tau^{(\Delta)}_p} V(q_{\tau^{(\Delta)}_p})=0$. Therefore, \eqref{verify3} implies 
	\begin{equation}
	\begin{split}\label{verify4}
	&\E \Big[ \int_0^{\tau^{(\Delta)}_p} e^{-rs} \Delta_s ds   \, \Big| \, \theta=1  \Big] 
	\leq V(q_0) =  \E \Big[ \int_0^{\tau_p^{(\alpha)}} e^{-rs} \alpha(q_s^{(\alpha)}) ds  \, \Big| \, \theta=1 \Big].
	\end{split}
	\end{equation}
	Finally, we conclude the optimality of $\alpha$ in \eqref{attacker_max} by \eqref{modified form attacker} and \eqref{verify4}. 
\end{proof}


The following proposition shows that in equilibrium, $S$ (the set of suspicion levels at which the defender stops the game) should have the form of $S=[p,1]$ for a constant $p>0$. 

\begin{proposition}\label{S_unique_appendix}
	Suppose that $((q_t)_{t\geq0},S,\alpha)$ is a Markov equilibrium. Then there exists a constant $p\in (0,1]$ such that $S=[p,1]$.
\end{proposition}
\begin{proof}
	We can easily see that if $q_0=0$ ($q_0=1$), then $\tau\equiv \infty$ ($\tau\equiv 0$) is the defender's optimal stopping time. This implies that $0\notin S$ and $1\in S$. Due to this observation and the closedness of $S$, to prove the proposition, it is enough to show that the set $S$ is connected. We prove it by contradiction. Suppose that there exist constants $0<\up<\op<1$ such that $\up,\op \in S$ and $(\up,\op)\cap S=\emptyset$.
	As in Subsection \ref{subsec_derive}, we derive the differential equation and the variational inequality for the value functions of the attacker and defender:
	\begin{align}
	&\textrm{If  }\,\, V'(q) \geq  -\tfrac{\sigma^2}{M(1-q)q}, \textrm{  then  } \begin{cases} 
	\frac{V''(q)}{2}+ \frac{V'(q)}{q} - \frac{r\sigma^2 V(q)}{M^2 q^2 (1-q)^2} + \frac{\sigma^2}{M q^2 (1-q)^2}=0,\\
	\alpha(q)=M, \\
	-r U(q)  +   \tfrac{q^2(1-q)^2 \alpha(q)^2 U''(q)}{2\sigma^2} + q \, \alpha(q)=0.
	\end{cases}   \label{character_DE1}\\
	&\textrm{If  }\,\, V'(q) <  -\tfrac{\sigma^2}{M(1-q)q}, \textrm{  then  } \begin{cases} 
	\frac{V''(q)}{2}- \frac{V'(q)}{1-q} - \frac{r}{\sigma^2} V'(q)^2 \,V(q)=0,\\
	\alpha(q)= -\tfrac{\sigma^2}{q(1-q)V'(q)},\\
	-r U(q)  +   \tfrac{q^2(1-q)^2 \alpha(q)^2 U''(q)}{2\sigma^2} + q \, \alpha(q)=0.
	\end{cases}   \label{character_DE2}
	\end{align}
	with the boundary conditions
	\begin{align}
	&V(\up)=V(\op)=0,\label{character_bd1}\\
	&U(\up)=l(1-\up), \quad U(\op)=l(1-\op),\label{character_bd2}\\
	&U'(\up)=U'(\op)=-l,   \label{character_bd3}\\ 
	&U(q) \leq l(1-q)  \textrm{  for  } q\in (\up,\op),   \label{character_bd4}\\
	&U(q) \leq \tfrac{M}{r}q, \,\, V(q)\geq 0  \textrm{  for  } q\in (\up,\op),   \label{character_bd5}
	\end{align}
	where the condition \eqref{character_bd3} is from the smooth-fit condition, and \eqref{character_bd5} is from the form of the optimization problems in \defref{def-equilibrium}.
	
	Now our goal is to show that there is no solution to the above system \eqref{character_DE1}-\eqref{character_bd5}. Suppose that there exist $V,\alpha, U$ satisfying \eqref{character_DE1}-\eqref{character_bd5}. Then, the following four steps lead us to a contradiction. 
	
	\medskip
	
	{\bf Step 1:} The set $\{ q\in (\up,\op): \, V'(q) \geq  -\tfrac{\sigma^2}{M(1-q)q}  \}$ is a connected set.
	
	{\it proof of Step 1.} Suppose that the set is not connected. Then, there exist constants $p_1$ and $p_2$ such that $\up<p_1<p_2<\op$ and 
	\begin{align}\label{step1_1}
	V'(q) <  -\tfrac{\sigma^2}{M(1-q)q} \textrm{  for  }q\in (p_1, p_2), \quad V'(p_1)= -\tfrac{\sigma^2}{M(1-p_1)p_1}, \quad V'(p_2)=  -\tfrac{\sigma^2}{M(1-p_2)p_2}.
	\end{align}
	The equalities in \eqref{step1_1} and the differential equation for $V$ in \eqref{character_DE2} produce
	\begin{align}\label{step1_2}
	V''(p_1)=-\tfrac{2\sigma^2 (M p_1 - r V(p_1))}{M^2 (1-p_1)^2 p_1^2},\quad V''(p_2)=-\tfrac{2\sigma^2 (M p_2 - r V(p_2))}{M^2 (1-p_2)^2 p_2^2}.
	\end{align}
	The equalities and inequality in \eqref{step1_1} imply that 
	\begin{align}\label{step1_3}
	\tfrac{d}{dq} \big( V'(q)+\tfrac{\sigma^2}{M(1-q)q} \big)\big|_{q=p_1}\leq 0, \quad 
	\tfrac{d}{dq} \big( V'(q)+\tfrac{\sigma^2}{M(1-q)q} \big)\big|_{q=p_2}\geq 0.
	\end{align}
	Combining \eqref{step1_2} and \eqref{step1_3}, we obtain $-\tfrac{\sigma^2 (M-2r V(p_1))}{M^2(1-p_1)^2 p_1^2} \leq 0$ and  $-\tfrac{\sigma^2 (M-2r V(p_2))}{M^2(1-p_2)^2 p_2^2} \geq 0$. These inequalities imply $V(p_1)\geq \tfrac{M}{2r}\geq V(p_2)$, but this contradicts to $V'(q)<-\tfrac{\sigma^2}{M(1-q)q}<0$ for $q\in (p_1,p_2)$. 
	
	\medskip

	{\bf Step 2:} There exists a constant $p^*\in (\up,\op)$ such that 
	\begin{align}
	\begin{cases}
	V'(q) \geq  -\tfrac{\sigma^2}{M(1-q)q}, &q\in [\up,p^*]\\
	V'(q) <  -\tfrac{\sigma^2}{M(1-q)q}, &q\in (p^*,\op]
	\end{cases}.
	\end{align}
	{\it proof of Step 2.} We first show that $V'(q) \geq  -\tfrac{\sigma^2}{M(1-q)q}$ for $q \geq \up$ close enough to $\up$. Suppose not. Then, Step 1 implies that $V'(q) <  -\tfrac{\sigma^2}{M(1-q)q}$ for $q> \up$ close enough to $\up$. The solution of the differential equation for $V$ in \eqref{character_DE2} with the boundary condition $V(\up)=0$ in \eqref{character_bd1} is
	$V(q)=\tfrac{\sigma}{\sqrt{r}} \varphi^{-1}(c_1 \cdot \tfrac{q-\up}{1-q})$ for a constant $c_1$. Then we reach a contradiction: $V\geq 0$ implies $c_1 \geq 0$, but 
	$$\tfrac{c_1 (1-\up)\sigma \sqrt{r}}{2\sqrt{r} (1-q)^2}e^{\varphi^{-1}(c_1 \cdot \frac{q-\up}{1-q})^2}=V'(q)< -\tfrac{\sigma^2}{M(1-q)q}<0
	$$
	implies $c_1<0$. Therefore, we conclude that $V'(q) \geq  -\tfrac{\sigma^2}{M(1-q)q}$ for $q \geq \up$ close enough to $\up$
	
	Suppose that $V'(q) \geq  -\tfrac{\sigma^2}{M(1-q)q}$ for $q\in [\up,\op]$. Then the differential equation for $U$ in \eqref{character_DE1} and the inequality \eqref{character_bd5} imply that $U''(q)\leq 0$ for $q\in [\up,\op]$. Since $U$ is concave, \eqref{character_bd2} and \eqref{character_bd4} imply that $U(q)=l(1-q)$ for $q\in [\up,\op]$, which does not satisfies the differential equation for $U$ in \eqref{character_DE1}. Therefore, we conclude that there exists $p^*\in (\up,\op)$ such that $V'(q) \geq  -\tfrac{\sigma^2}{M(1-q)q}$ for $q\in [\up,p^*]$, and $V'(q) <  -\tfrac{\sigma^2}{M(1-q)q}$ for $q>p^*$ close enough to $p^*$.
	Furthermore, Step 1 implies that $V'(q)  <  -\tfrac{\sigma^2}{M(1-q)q}$ for $q\in (p^*,\up]$.

	\medskip
	
	{\bf Step 3:} $U''(q)\leq 0$ for $q\in (\up,p^*)$ and $U''(q)<0$ for $q\in [p^*,\op)$.\\
	{\it proof of Step 3.} We first observe that \eqref{character_bd2} and \eqref{character_bd5} imply 
	\begin{align}\label{up_ineq}
	r l (1-p^*)- M p^*< r l (1-\up)- M \up \leq 0.
	\end{align}
	We rewrite the differential equation for $U$ in \eqref{character_DE1} and \eqref{character_DE2} as
	\begin{align}
	U''(q)=\tfrac{2\sigma^2 (r U(q)-q \alpha(q))}{(1-q)^2 q^2 \alpha(q)^2}.
	\end{align}
	Therefore, it is enough to check that $r U(q)-q \alpha(q) \leq 0 $ for $q\in (\up,p^*)$ and $r U(q)-q \alpha(q) < 0 $ for $q\in [p^*,\op)$. \\
	(i) For $q\in (\up, p^*)$, by the result in Step 2 and \eqref{character_DE1} and \eqref{character_bd5}, we have $r U(q)-q \alpha(q)\leq 0$.\\
	(ii) For $q\in [p^*,\op)$, by the result in Step 2 and \eqref{character_DE2} and \eqref{character_bd2}, we have $r U(q)-q \alpha(q)\leq r l(1-q)-q \alpha(q)$. To prove $r l(1-q)-q \alpha(q)< 0$ for $q\in [p^*,\op)$, it is enough to show that the function $\tfrac{r l(1-q)-q \alpha(q)}{1-q}$ decreases in $q$ on $(p^*,\op)$, because $r l(1-p^*)-p^* \alpha(p^*)< 0$ by \eqref{up_ineq}. Indeed, for $q\in (p^*,\op)$, the solution of the differential equation for $V$ in \eqref{character_DE2} with \eqref{character_bd1} is the form of $V(q)=\tfrac{\sigma \varphi^{-1}(c_2 \cdot \frac{\op-q}{1-q})}{\sqrt{r}}$ for a constant $c_2>0$, so using this, we obtain
	$$
	\big(\tfrac{r l(1-q)-q \alpha(q)}{1-q}\big)'= -\tfrac{2\sqrt{r}\sigma \varphi^{-1}(c_2 \cdot \frac{\op-q}{1-q})}{(1-q)^2}<0 \quad \textrm{for  }q\in (p^*,\op).
	$$
	{\bf Step 4:} The system \eqref{character_DE1}-\eqref{character_bd5} does not have a solution.\\
	{\it proof of Step 4.} In Step 3, we concluded that if there exists a solution to the system \eqref{character_DE1}-\eqref{character_bd5}, then $U$ should be a concave function on $[\up,\op]$. The conditions \eqref{character_bd2} and \eqref{character_bd4}, together with the concavity of $U$, imply that $U(q)=l(1-q)$ for $q\in [\up,\op]$. This contradicts the result in Step 3: $U''(q)<0$ for $q\in [p^*,\op)$.
\end{proof}

\section{Estimation of the Initial Suspicion Level $q_0$}
In this section, we present empirical and theoretical estimations regarding the initial suspicion level $q_0$, interpreting it as the fraction of infected users among a whole population of users whose traffics are observed by the defender. Estimation of the initial suspicion level is essential in several industries related to cyber-security. For instance, an actuary may have to estimate how many users will be infected in the future when designing a cyber-insurance contract. However, estimation of $q_0$ using actual data may not be reliable due to insufficient amount of actuarial data \cite{WanchunDou}. In this case, our proposed methods can complement the estimation.

\subsection{Empirical Estimation of $q_0$}
In this subsection, we propose a method of empirically estimating the actual fraction of infected users using data that contain whether a user is blocked by the defender before random termination. 

We assume that the defender's perception $q_0$ might be different from the actual fraction $x$. To be specific, we consider a situation that 
\begin{align}
\PP_{defender}(\theta=1)=q_0 \neq \PP_{true}(\theta=1)=x,
\end{align}
where $\PP_{defender}$ is the probability measure describing the defender's belief and $\PP_{true}$ is the actual probability measure. The proposition below provides an unbiased estimator of $\PP_{true}(\theta=1)=x$ in terms of the ratio of the blocked suspects when the defender has the belief $\PP_{defender}(\theta=1)=q_0$.

\begin{proposition}\label{prop_7_simple}
Suppose that there are $N\in \mathbb{N}$ suspects. For suspect $i\in \{ 1,2,\ldots, N\}$, $\theta^{(i)}$ indicates whether suspect $i$ is an attacker or not, $(W_t^{(i)})_{t\geq 0}$ represents the noise in \eqref{Y}, and $T^{(i)}$ is the random termination time in \eqref{T_exp}. We assume that these random variables and Brownian motions are all independent under $\PP_{defender}$ and $ \PP_{true}$, and
\begin{align}
\PP_{true}(\theta^{(i)}=1)=x \quad \textrm{and} \quad \PP_{defender}(\theta^{(i)}=1)=q_0\in (0,p), \quad 1\leq i \leq N.
\end{align}
Each suspect and the defender play the game under the probability measure $\PP_{defender}$, and we denote by $\tau_p^{(i)}$ the first hitting time for suspect $i$. Then, $\mu_\theta$ defined below is an unbiased estimator of $x$:
\begin{equation}\label{eqn-x}
\mu_\theta := \frac{p(1-q_0)}{(p-q_0)u(q_0)} \cdot \frac{1}{N}\sum_{i=1}^N 1_{\left\{ \tau_p^{(i)}<T^{(i)} \right\}}-\frac{q_0(1-p)}{p-q_0},
\end{equation}
where
\begin{equation}\begin{split}\label{block_prob}
u(q_0)=\begin{cases}
\big(\frac{q_0(1-p)}{p(1-q_0)}\big)^a, & \textrm{if  }\,\,\frac{r \sigma^2}{M^2} \geq 1\\
\frac{acM\sqrt{\pi}}{2\sigma\sqrt{r}}\big(\frac{q_0(1-q^*)}{q^*(1-q_0)}\big)^a, & \textrm{if  }\,\, \frac{r \sigma^2}{M^2} < 1 \,\, \& \,\, q_0\in (0,q^*] \\
\frac{p(1-q_0)}{q_0(1-p)} e^{-y(q_0)^2} - c \sqrt{\pi} \,y(q_0), & \textrm{if  }\,\, \frac{r \sigma^2}{M^2} < 1\,\, \& \,\, q_0\in (q^*,p)
\end{cases}
\end{split}
\end{equation}
and the constants $a,c,$ and $q^*$ and the function $y$ are defined in \eqref{const}.
\end{proposition}

\begin{proof}
By straightforward computations, we check that $u$ in \eqref{block_prob} satisfies $u\in C^2([0,p))$ and solves the following differential equation:
\begin{equation}\label{u_ode}
\begin{cases}-r u(q)+\tfrac{q^2(1-q)^2 \alpha(q)^2}{2\sigma^2} u''(q) +\tfrac{q(1-q)^2 \alpha(q)^2}{\sigma^2} u'(q)=0\\
u(0)=0,\quad u(p)=1.
\end{cases}
\end{equation} 
Conditioned on $\theta^{(i)}=1$ and the initial suspicion level $q_0$, Ito's formula and \eqref{u_ode} produce the following:
\begin{equation}\label{u_ito}
e^{-r(t\wedge \tau_p^{(i)})} u(q_{t\wedge \tau_p^{(i)}})= u(q_0) + \int_0^{t \wedge \tau_p^{(i)}} e^{-rs} \tfrac{q_s(1-q_s)\alpha(q_s)}{\sigma} u'(q_s) dW_s^{(i)}.
\end{equation} 
In \eqref{u_ito}, the stochastic integral is a square integrable martingale since the integrand is bounded. Therefore,
\begin{equation}\label{u_ito2}
\begin{split}
u(q_0)&=\lim_{t\to \infty}\E_{defender}\left[e^{-r(t\wedge \tau_p^{(i)})} u(q_{t\wedge \tau_p^{(i)}})\, \big| \,  \theta^{(i)}=1\right]\\
&=\E_{defender}\left[e^{-r \tau_p^{(i)}} \cdot 1_{\{\tau_p^{(i)}<\infty\}}\, \big| \,  \theta^{(i)}=1\right]\\
&=\PP_{defender}\left(\tau_p^{(i)}<T^{(i)} \, \big| \, \theta^{(i)}=1\right),
\end{split}
\end{equation} 
where the second equality is due to the dominated convergence theorem and $u(p)=1$, and the third equality holds since $T^{(i)}$ is exponentially distributed and independent of other random variables. 
By the same way, we also obtain the expression of $\PP_{defender}(\tau_p^{(i)}<T^{(i)} \, |\,  \theta=0)$, and the result is summarized below:
\begin{equation}\begin{split}\label{u_expss}
&\PP_{defender}(\tau_p^{(i)}<T^{(i)} \, |\, \theta=1) = u(q_0),\\
&\PP_{defender}(\tau_p^{(i)}<T^{(i)} \, | \, \theta=0) =\tfrac{q_0(1-p)}{p(1-q_0)} \, u(q_0).\\
\end{split}
\end{equation}

Since $\PP_{true}(\theta^{(i)}=1)=x$, the probability of user $i$ being blocked before the random termination is
\begin{equation}\label{eqn-prob-taulessT}
\begin{split}
\PP_{true}(\tau_p^{(i)} < T^{(i)} ) &= \PP_{true}(\tau_p^{(i)}<T^{(i)} \, | \, \theta^{(i)} = 1) \cdot x  + \PP_{true}(\tau_p^{(i)}<T^{(i)} \, | \, \theta^{(i)} = 0) \cdot (1 - x)\\
&= \PP_{defender}(\tau_p^{(i)}<T^{(i)} \, | \, \theta^{(i)} = 1) \cdot x  + \PP_{defender}(\tau_p^{(i)}<T^{(i)} \, | \, \theta^{(i)} = 0) \cdot (1 - x),
\end{split}
\end{equation}
where the second equality is due to the observation that $\PP_{true}=\PP_{defender}$ once the value of $\theta^{(i)}$ is realized. 

Finally, we combine \eqref{u_expss} and \eqref{eqn-prob-taulessT} to conclude that $\EE_{true}[\mu_\theta]=x$. 
\end{proof} 
 
The empirical estimation method presented in \proref{prop_7_simple} does not require inspection results. That is, the data set does not need to indicate whether a blocked suspect is actually an attacker or an innocent user. This method is useful in the sense that a defender can adjust the probability of cyber-attacks before inspection results come out.

We can actually relax the assumption that all suspects and the defender share the same suspicion level $q_0$ at the beginning of the game. Even if suspect $i$ is an attacker with probability $q_0^{(i)}$, one can still utilize the expression in \eqref{p alpha form}. The defender chooses a certain $q_0$, and construct observation data starting from $q_0$. In this case, the estimator $\mu_\theta$ in Proposition \ref{prop_7_simple} is the ratio of infected users to the whole population of users.

\subsection{Theoretical Estimation of $q_0$}

We move onto a theoretical estimation on the probability of cyber-attacks. So far, we have assumed that the initial suspicion level $q_0$, which can be interpreted as the probability of cyber-attacks, is taken as given. However, we can imagine cases in which the attacker chooses the attack probability $q_0$ as in \cite{LinChen2009,HaoWu}. With the probability $q_0$, the attacker actually launches cyber-attacks, whereas with the complementary probability $1-q_0$, the attacker leaves the game. For example, one may imagine a bot herder, a malicious hacker who controls a botnet (many bot-infected devices) to attack a target and chooses the proportion of active bots. In this case, $q_0$ represents the proportion of bots that actually launch attacks.

Suppose that the attacker chooses to attack with the probability $q_0 \in [0,1]$.  When the attacker actually attacks, her expected profit $V(q_0)$ in equilibrium is
\[
V(q_0) = \mathbb{E}\left[ \int_0^{T\wedge \tau_p}  \alpha(q_t) dt  \Big\vert \theta = 1  \right].
\]
Then, $q_0 V(q_0)$ is the attacker's expected profit when she chooses to attack with the probability $q_0$. Therefore, the \emph{optimal attack probability} $\hat{q}$,
\begin{align}\label{hq_def}
\hat q \in \argmax_{q_0 \in [0,1]} q_0 V(q_0),
\end{align}
can be a reasonable theoretical estimation of $q_0$ when we interpret $q_0$ as the attack probability. The proposition below characterizes the optimal attack probability.

\begin{proposition}\label{prop-critical-l}
There exists a unique optimal attack probability $\hat q\in (0,p)$ in \eqref{hq_def}. 
\end{proposition}
\begin{proof}
We first observe the sign of $\tfrac{\partial^2}{\partial q^2} \big( q V(q) \big)$, whose expression is
\begin{align}\label{qV''}
\begin{cases}
-\tfrac{a(1+a)}{(1-q)^2 q}\left(\tfrac{M}{r}\right) (\tfrac{1-p}{p})^a (\tfrac{q}{1-q})^a,  &\textrm{   if     } \tfrac{r \sigma^2}{M^2} \geq 1 \textrm{  and  } q\in [0,p)\\
-\tfrac{(1+a)\sigma^2}{M(1-q)^2 q}(\tfrac{1-q^*}{q^*})^a (\tfrac{q}{1-q})^a,  &\textrm{   if     } \tfrac{r \sigma^2}{M^2} <1 \textrm{  and  } q\in [0,q^*]\\
-\tfrac{\sigma(1-p)c\sqrt{\pi}e^{2y(q)^2}\left(2p(1-q)e^{-y(q)^2}-c \sqrt{\pi}(1-p)q y(q)  \right)}{2p^2 (1-q)^4\sqrt{r}},  &\textrm{   if     } \tfrac{r \sigma^2}{M^2} <1 \textrm{  and  } q\in (q^*,p)
\end{cases}.
\end{align}
The inequalities in \lemref{tech_lemma} ensure that the above three expressions are all strictly negative. Therefore, the map $q\mapsto q V(q)$ on $q\in [0,p]$ is strictly concave, and we conclude that there exists unique $\hq$ such that
$$\hat{q} = \argmax_{q \in [0,p]} q V(q).$$
Since $q V(q)=0$ for $q\in [p,1]$ and $q V(q)>0$ for $q\in (0,p)$, we conclude that $\hq$ above also satisfies \eqref{hq_def}.
\end{proof}



In many studies, expected losses due to cyber-attacks are taken as given \cite{WanchunDou,HaoWu}. However, the expected losses are related to the probability of cyber-attacks in dynamic environments. This is because if the probability of cyber-attacks is lower, the attacker can gain more due to longer duration of the game. Our model suggests a relationship between the attack probability and the expected losses. Given the optimal attack probability $\hat{q}$, the expected losses would be $V(\hat{q})$.




\section{Calibration of Model Parameters}
We calibrate the model parameters based on the report conducted by Ponemon Institute and sponsored and published in 2020 by IBM Security.\footnote{https://www.ibm.com/security/digital-assets/cost-data-breach-report. This report is referred to as ``the report" throughout this section.} The report collected the information on data breach from more than 500 organizations. To name a few items in the report, there are components of costs, detection times by industry and by nation, and root causes of data breaches. In the report, a main reason for data breach falls into one of the three categories: system glitches, human errors, and malicious attacks. Among these categories, 50\% of data breaches are due to malicious attacks, 13\% of which are carried out by nation state attackers. Even though our model is inspired by APTs, our model is intended to capture long-term data breach including some characteristics of APTs. For the calibration purpose, we use the aggregate data as the report does not provide detailed information on single data breaches. Nevertheless, we believe that the report would fit into our model to some extent.

Before we describe our calibration, it would be worth discussing an issue on discrepancy between our model and the report. For instance, let us consider security automation deployment in the report. It refers to enabling augment or replace human intervention in the identification and containment of cyber exploits or breaches. According to the report, average security automation deployment by industry ranges from 49\% to 68\%. This means that a certain fraction of industry does not deploy security automation at all. No deployment of security automation can be a reason for a longer average detection time because firms without security automation deployment possibly detect data breaches at later times. The report does not provide details of data composition and the defender in our model can be thought of as a security automation system, which means that our calibration is likely to result in inaccurate values to some extent. It would be a future project to develop a model that can handle the discrepancy issues.\footnote{One might be wondering why the report is chosen. A main reason is because the report deals with data breaches in the long run as our model deals with persistent data leakage over time. Another good candidate set of data for calibration would be the PRC database. Although the PRC database provides information on types of data breaches, its composition may not be perfectly suited to our setting. According to the PRC database, the most frequent type of data breach is ``PHYS", paper documents that are lost or stolen, accounting for almost 1,400 incidents. The second most frequent type of data breach is ``DISC", unintended disclosure that is not involved with hacking and intentional breaches, accounting for slightly more than 1,000 incidents. ``HACK" type, which means `hacked by an outside party or infected by malware', is the third most frequent type of data breach and accounts for around 900 incidents. It seems that the two most frequent types  PHYS and DISC  are possibly related to one-shot events.
}

There are five model parameters to calibrate: $l, M, q_0, r$ and $\sigma$. We substitute the false alarm cost $l$ with the lost business cost of \$1.52M in the report. Lost business includes business disruption, opportunity cost during system downtimes, lost costumers, and reputation losses. These events can happen when the defender falsely shuts down servers, and thus the lost business cost is used as the false alarm cost. For the upper bound of the attack intensity $M$, we set a sufficiently large number $M=100$ as it has little impacts on equilibrium.\footnote{See the graph of $\alpha$ for varying $M$ in Figure \ref{alpha_fig} and the graph of $p$ as a function of $M$ in Figure \ref{p_fig}.} For the initial attack probability $q_0$, we use the optimal attack probability $\hat q$ in equation \eqref{hq_def}, which is a function of the remaining two undetermined parameters $r$ and $\sigma$.

Among information contained in the report, we focus on the {\it average cost of a data breach} and {\it average detection time}, by industry. We set $V(\hat q)$ as the average cost of a data breach, and $\mathbb{E}[\tau_p \wedge T \vert \theta = 1, q_0=\hat q]$ as the average detection time. The expression of the function $V$ is given in \proref{solutions}, and the expression in equation \eqref{u_ito2} provides
\begin{align}
\mathbb{E}[\tau_p \wedge T \vert \theta = 1,\, q_0=\hat q] &= \E \left[ \int_{0}^{\tau_p} e^{-rt} dt  \Big \vert \theta = 1,\, q_0=\hat q  \right]= \frac{1-\E[e^{-r\tau_p} \vert \theta = 1,\, q_0=\hat q]}{r}\nonumber\\
&= \frac{1-u(\hat q)}{r},\nonumber
\end{align}
where the function $u$ is defined in equation \eqref{block_prob}. We calibrate $r$ and $\sigma$ to match the values of $V(\hat q)$ and $\mathbb{E}[\tau_p \wedge T \vert \theta = 1, q_0=\hat q]$.

\begin{table}[t]
	\centering
	\begin{tabular}{ | c | c |  c |  c  |  c |  c |}
		\hline
		Industry & Average cost & Average detection time  & $r$ &  $\sigma$ &  $\alpha(\hat q)/\sigma$ \\
		\hline
		Healthcare       & \$7.13M   &  329 days  &  0.32   &  7.1  &   1.09  \\
		Energy             & \$6.39M   &  254 days  &  0.42   &  7.3  &   1.26  \\
		Financial          & \$5.85M   &  233 days  &  0.46   &  6.9  &   1.32  \\
		Pharmaceuticals & \$5.06M &  257 days  &  0.42   &  5.7  &   1.26 \\
		Technology      & \$5.04M   &  246 days  &  0.44   &  5.8   &  1.29 \\
		Manufacturing          & \$4.99M   &  302 days  & 0.35    &  5.1   &  1.16  \\
		Services           & \$4.23M   &  286 days  & 0.38    &  4.4  &  1.21 \\
		Entertainment   & \$4.08M   &  314 days  & 0.35   &  4.1  &   1.16 \\
		Education         & \$3.90M   &  283 days  & 0.39   &  4.1   &   1.23 \\
		Transportation  & \$3.58M   &  275 days  & 0.40   &  3.8  &  1.25 \\
		Communication & \$3.01M   & 251 days  & 0.44   &  3.3   &  1.33 \\
		Consumer         & \$2.59M   & 307 days  & 0.36    & 2.5   &  1.21 \\
		Retail                & \$2.01M   & 311 days   & 0.37    & 1.9   &  1.25 \\
		Hospitality        & \$1.72M   & 275 days   & 0.43    & 1.7   &  1.36 \\
		Media               & \$1.65M   & 281 days   & 0.42    & 1.6   &  1.35 \\
		Research         & \$1.53M   & 244 days   & 0.49    & 1.6   &  1.47 \\
		Public               & \$1.08M   & 324 days   & 0.39    & 0.9   &  1.36 \\
		\hline
		Global average & \$3.86M   & 280 days   & 0.39    & 4.1  &  1.23 \\
		\hline
	\end{tabular}
	\caption{Calibration of Model Parameters for 17 industries}\label{table-byindustry}
\end{table}

Table \ref{table-byindustry} presents the calibrated parameters of the model for 17 industries used in the report and provides basic building blocks for discussions about the calibration of the model parameters.

For discussion of calibrated $r$, we would like to extend the meaning of $r$. As mentioned in Section 2, $r$ represents how often random termination occurs. However, it is worth mentioning that $r$ can represent the overall time preference. That is, $r$ can also represent how fast values of data decay over time or foregone values of other alternatives. For instance, a large value of $r$ can imply a high frequency of detecting cyber-attacks (or equivalently, a high risk of the attacker being detected), a fast decay of current data value, and a large profits of the attacker doing other things rather than cyber-attacks. For this reason, we refer to $r$ as \emph{time preference factor} in this section.

According to our calibration, two industries with the highest time preference factors are financial industry and research industry. Intuitively, one can think of high time preference factors in these industries. For research industry, frontier researches and state-of-the-art technologies are rapidly devalued over time. For financial industry, personally identifiable information of important customers would depreciate quickly, compared to other industries. In addition, hackers would take larger risks in finance and research industries as information and data are protected by better systems in these industries.

Three industries with the lowest time preference factors are healthcare industry, entertainment industry, and manufacturing industry. Information on patients in healthcare industry, consumers in manufacturing and entertainment industries is important, and values of information may not change quickly over time. Therefore, one can anticipate low time preference factors in those industries.

We close this section by discussing the noise intensity $\sigma$ and the informativeness $\frac{\alpha(\hat q)}{\sigma}$ of the signal process. The noise intensity represents the volatility of normal activities, and Table \ref{table-byindustry} shows that industry with the lowest noise intensity is the public sector, as one may expect. Financial, healthcare, and energy industries have the highest noise intensity. This would be because financial, healthcare, and energy services are necessary in daily lives and thus there are many transactions and emergency cases in those industries. For $\frac{\alpha(\hat q)}{\sigma}$, we find that it is rather stable across all industries. This would be a model property as we conjecture that there is an optimal range of $\frac{\alpha}{\sigma}$ for the attacker.

\section{Numerical Illustrations}

In this section, we graphically illustrate the results in the previous sections for varying parameters. 
Based on the discussions in Section 6 and the ``Global average" row in Table \ref{table-byindustry}, we set the model parameters as
\begin{align}
M=100, \quad l=1.52, \quad r=0.39, \quad \sigma=4.1.
\end{align}
For these parameters, the equilibrium threshold is $p=0.29$ and the optimal attack probability is $\hat q=0.14$. 

\begin{figure}[h]
	\small{	\begin{center}$
			\begin{array}{cc}
			\includegraphics[width=0.4\textwidth]{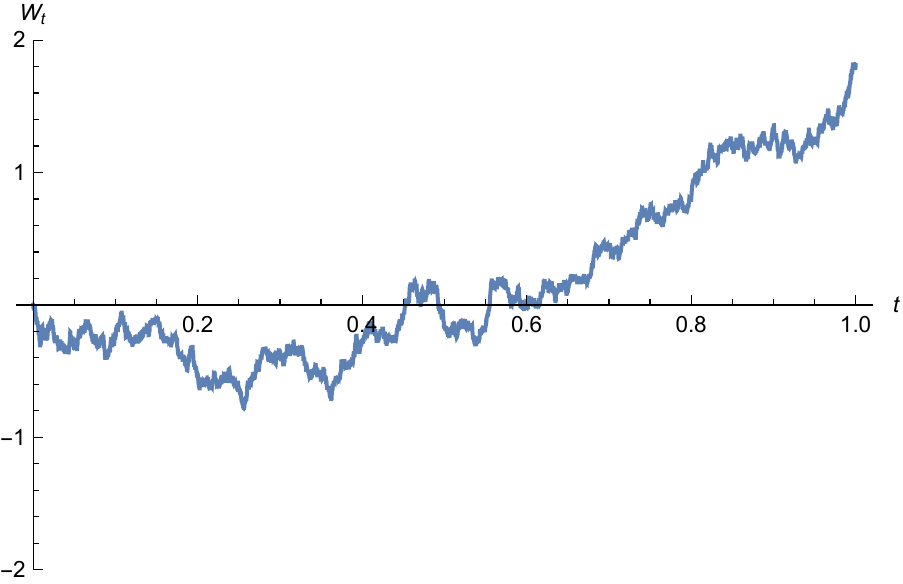} &
			\includegraphics[width=0.4\textwidth]{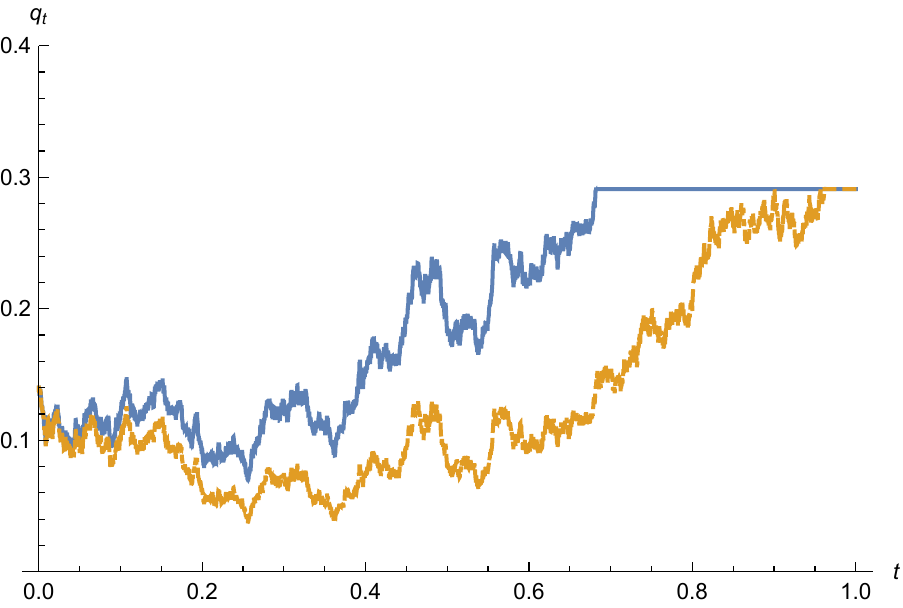} 		
			\end{array}$
	\end{center}}
	\vspace*{-5mm}
	\caption{The left graph is a sample path of a standard Brownian motion. The right graph is the corresponding paths of the suspicion level process when $\theta=1$ (solid line) and $\theta=0$, with parameters
$q_0=0.14, \, M=100, \, \sigma=4.1,\, r=0.39,\, l=1.52$. In this case, the equilibrium stopping threshold $p$ is $0.29$.	}
	\label{block_path}
\end{figure}

The left graph in Figure \ref{block_path} presents a sample path of the noise $(W_t)_{t\geq 0}$, and the right graph presents the corresponding paths of the suspicion level process $(q_t)_{t\geq 0}$, when $\theta=1$  (solid line) and $\theta=0$ (dashed line), respectively. As the definition of $q_t$ and the SDE \eqref{eqn-suspicion-level} indicate, for a given sample path of $W$, we observe that $q_t$ for $\theta=1$ case is always higher than $q_t$ for $\theta=0$ case. Therefore, the defender blocks the suspect earlier ($q_t$ hits the stopping threshold $p$ earlier) when the suspect is the attacker, as the figure shows.

\begin{figure}[h]
	\small{	\begin{center}$
			\begin{array}{cc}
			\includegraphics[width=0.3\textwidth]{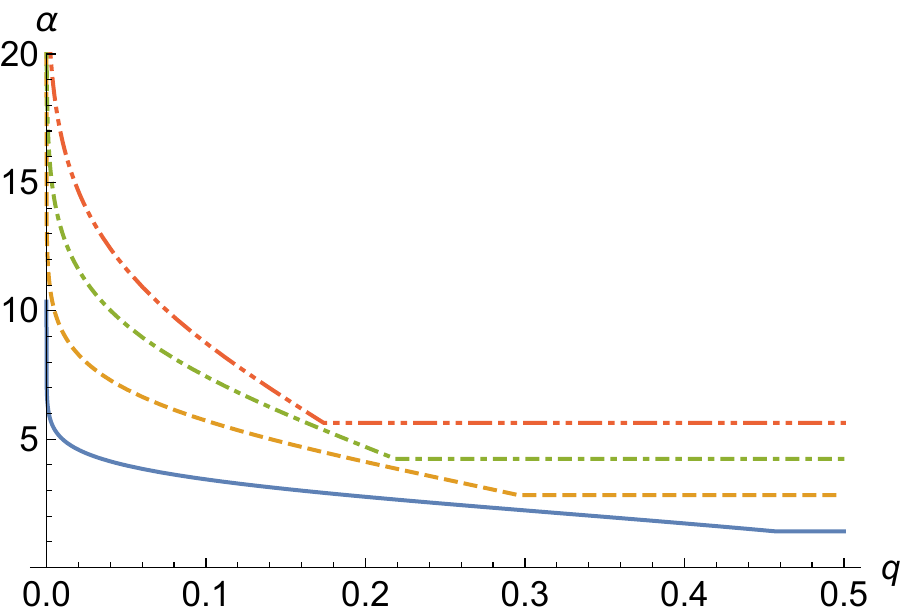} &
			\includegraphics[width=0.3\textwidth]{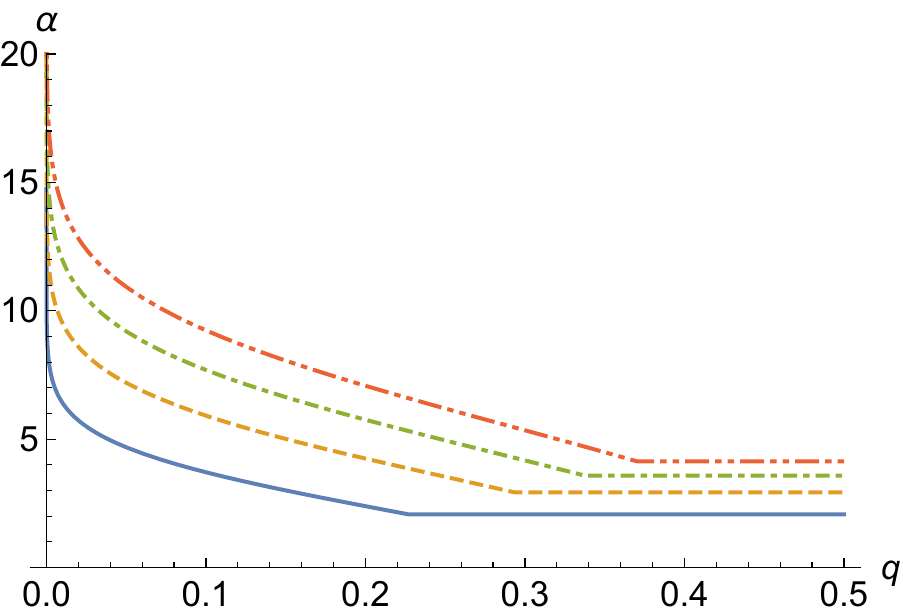}  \\
			\sigma=2 \, \text{(---)},\, \sigma=4 \, \text{(-\,-)},    &  r=0.2 \, \text{(---)},\, r=0.4 \, \text{(-\,-)},\\
			\sigma= 6  \, \text{(-$\cdot$-)},  \sigma=8,  \,  \text{(-$\cdot\cdot$-)},  & r=0.6\, \text{(-$\cdot$-)}, \, r=0.8 \, \text{(-$\cdot\cdot$-)},\\
			M=100, \, r=0.39,\, l=1.52 & M=100, \sigma=4.1, \, l=1.52 \\
			\includegraphics[width=0.3\textwidth]{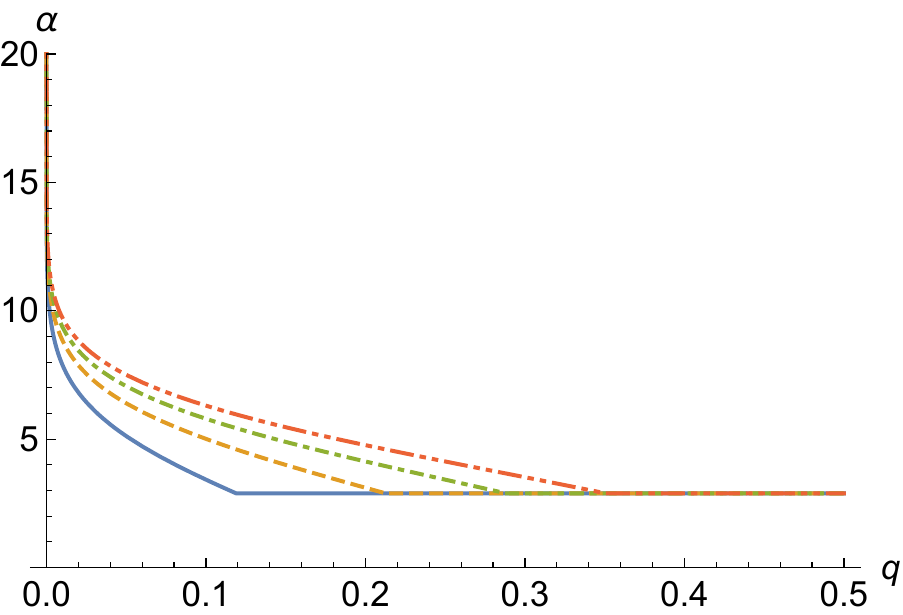} &
			\includegraphics[width=0.3\textwidth]{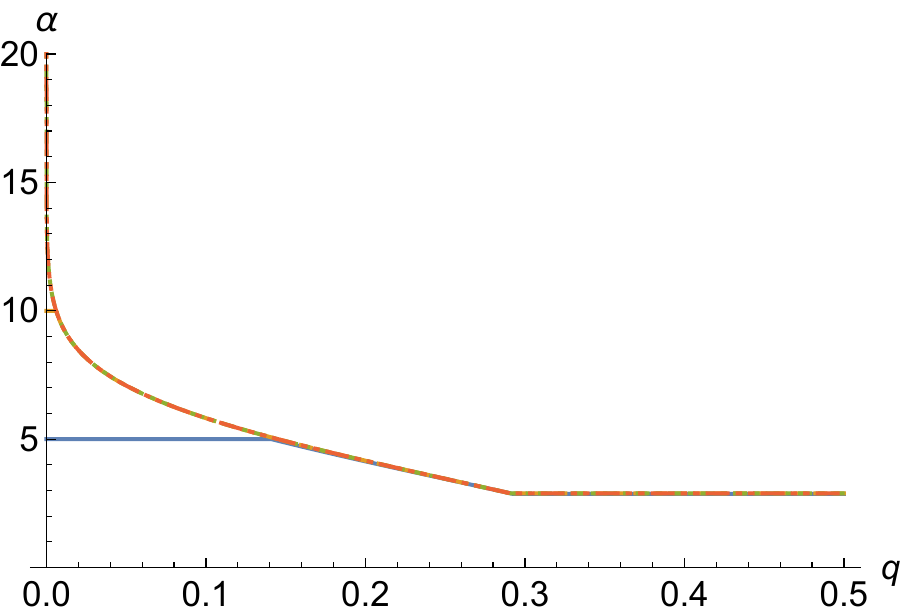} \\ 
			l=0.5 \, \text{(---)},\, l=1.0 \, \text{(-\,-)},   &  M=5\, \text{(---)},\, M=10  \, \text{(-\,-)},  \\
		        l=1.5 \,\text{(-$\cdot$-)}, \, l=2.0  \,  \text{(-$\cdot\cdot$-)}, & M=20 \, \text{(-$\cdot$-)}, \, M=100 \,  \text{(-$\cdot\cdot$-)},\\
			M=100, \, r=0.39,\, \sigma=4.1 & r=0.39, \sigma=4.1, \, l=1.52
			\end{array}$
	\end{center}}
	\vspace*{-5mm}
	\caption{Graphs of equilibrium attack intensities for varying $l,r,\sigma$, and $M$. } 
	\label{alpha_fig}
\end{figure}

\begin{figure}[h]
	\small{	\begin{center}$
			\begin{array}{cc}
			\includegraphics[width=0.3\textwidth]{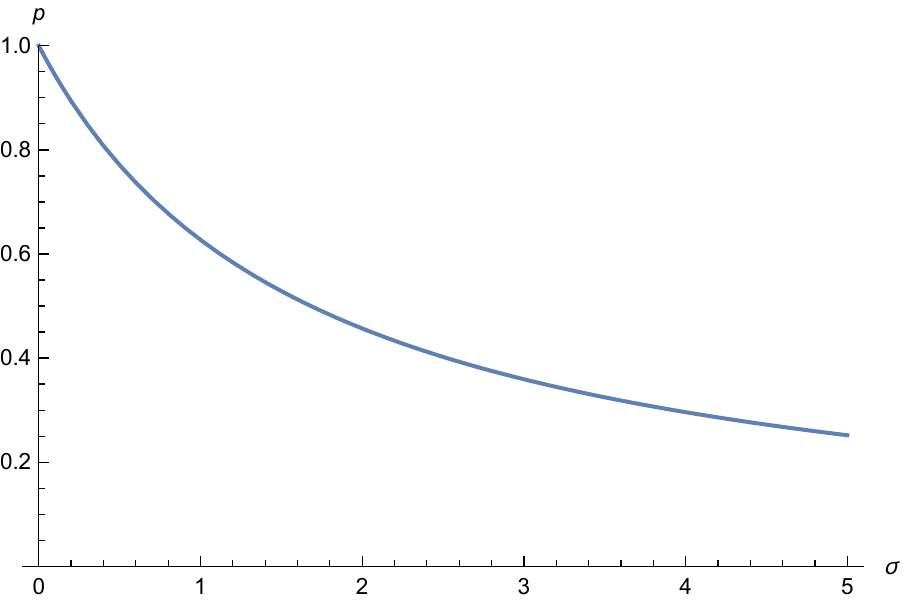} &
			\includegraphics[width=0.3\textwidth]{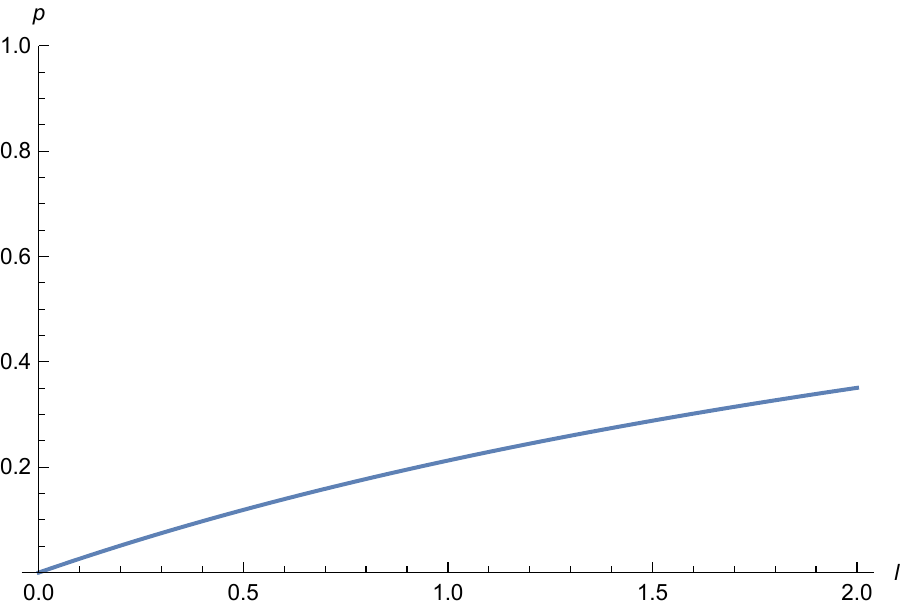} \\
			M=100,\, r=0.39,\, l=1.52 & M=100,\, r=0.39,\, \sigma=4.1 \\
			 \includegraphics[width=0.3\textwidth]{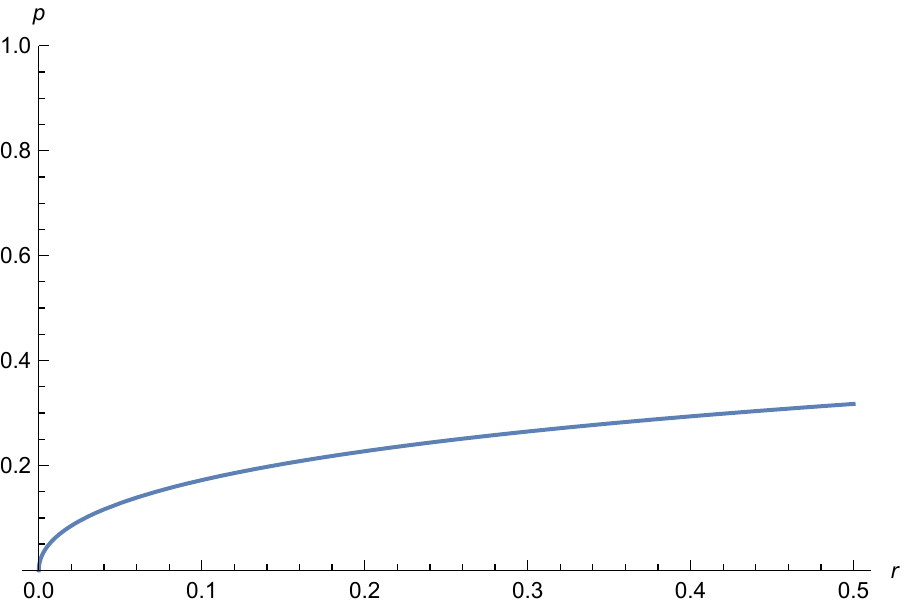} &
			\includegraphics[width=0.3\textwidth]{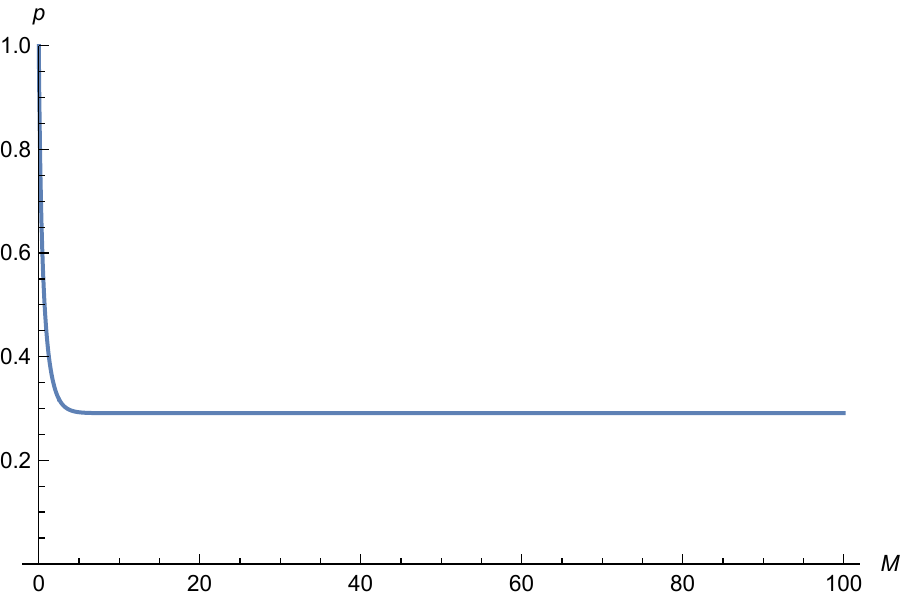} \\
			 M=100, \, \sigma=4.1,\, l=1.52 & r=0.39, \, \sigma=4.1,\, l=1.52\\	
			\end{array}$
	\end{center}}
	\vspace*{-5mm}
	\caption{The equilibrium threshold $p$ for varying $l,r,\sigma$, and $M$. } 
	\label{p_fig}
\end{figure}

The equilibrium attack intensity function $\alpha$ in \eqref{p alpha form} is illustrated in Figure \ref{alpha_fig}. One can observe in the figure that the attacker chooses the highest attack intensity $M$ when the suspicion level $q_t$ is low enough. The attack intensity function  decreases as the suspicion level $q_t$ increases. The intuition behind this is the following. As the suspicion level increases, the duration of the game decreases. The attacker decreases the attack intensity to slow down the increase of $q_t$ (see the form of the SDE \eqref{eqn-suspicion-level}) and induces the defender to stop the game later. This strategic behavior of our attacker can be interpreted as reducing the current profit to extend the duration of the game.

The comparative statics in \proref{coro-alpha-comp-statics} is demonstrated in Figures  \ref{alpha_fig} and \ref{p_fig}.  For the parameters we chosen, the graphs show that the equilibrium stopping threshold $p$ increases in $l$ and $r$ and decreases in $\sigma$ and $M$, and the equilibrium attack intensity  $\alpha$ increases in $l,r,\sigma$ and $M$.

\begin{figure}[h]
	\small{	\begin{center}$
			\begin{array}{cc}
			\includegraphics[width=0.4\textwidth]{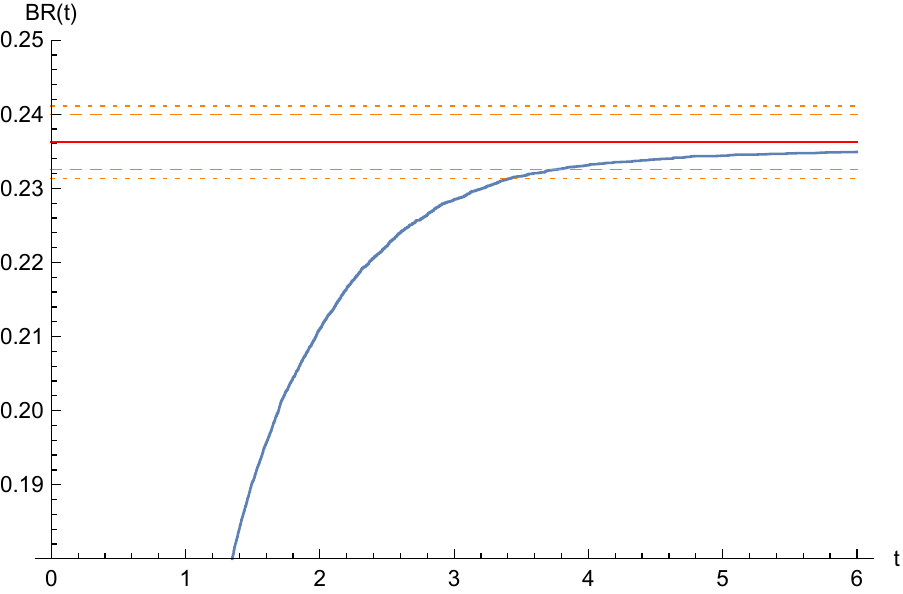} &
			\includegraphics[width=0.4\textwidth]{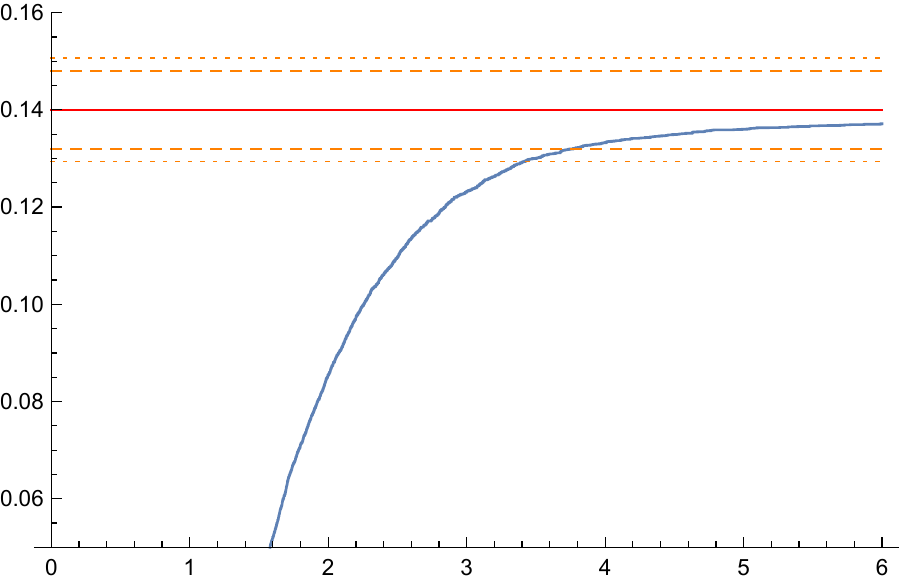} 		
			\end{array}$
	\end{center}}
	\vspace*{-5mm}
	\caption{The left graph is ratio of the blocked suspects as time goes (see \eqref{blocked_ratio} for $BR(t)$) with $N=50,000$ simulations where the solid red line is the desired ratio computed by \eqref{eqn-x}.  The right graph is the corresponding empirical estimation of $x$ (see \eqref{true_x} for $EE(t)$), where the solid red line is the true value of $x=0.14$. The parameters are $M=100, \, r=0.39, \sigma=4.1,\, l=1.52$, and the solution $q_t$ of \eqref{bayesian} is computed with the initial prior $q_0=0.1$. In each graph, the dashed line and the dotted line correspond to 95\% and 99\% confidence intervals, respectively. 
	}
	\label{simulation_fig}
\end{figure}

Figure \ref{simulation_fig} is a simulation result illustrating  the empirical estimator $\mu_\theta$ in \proref{prop_7_simple}. For given $x=0.14$, we generate
 $N=50,000$ sample paths of standard Brownian motion, $N$ realizations of the random termination time $T$ (from exponential distribution in \eqref{T_exp}), and $N$ realizations of the suspect types (from Bernoulli distribution with parameter $x=0.14$). For these $N$ realizations, we pick $q_0=0.1$ and calculate the suspicion level processes $q_t$ starting at $q_0=0.1$. If a suspicion level process ceases due to the random termination time or hits the equilibrium threshold $p$, then the suspicion level process is frozen. In these $N$ samples, we count the number of cases that $q_t$ hits $p$ before the random termination. Then we use \eqref{eqn-x} to infer the actual initial suspicion level. The left graph describes the ratio of the blocked suspects, as a function of $t$:
\begin{align}\label{blocked_ratio}
BR(t):=\frac{1}{N}\sum_{i=1}^N 1_{\left\{ \tau_p^{(i)}<T^{(i)} \wedge t \right\}},
\end{align}
The right graph is the corresponding empirical estimation, as a function of $t$:
  \begin{align}\label{true_x}
EE(t) := \frac{p(1-q_0)}{(p-q_0)u(q_0)} \cdot \frac{1}{N}\sum_{i=1}^N 1_{\left\{ \tau_p^{(i)}<T^{(i)}\wedge t \right\}}-\frac{q_0(1-p)}{p-q_0},
\end{align}
  In Figure \ref{simulation_fig}, we can see that the empirical estimation $EE(t)$ approaches to the true value $x=0.14$ as time goes on.

\begin{figure}[h]
	\begin{center}
			$\begin{array}{c}
			\includegraphics[width=0.4\textwidth]{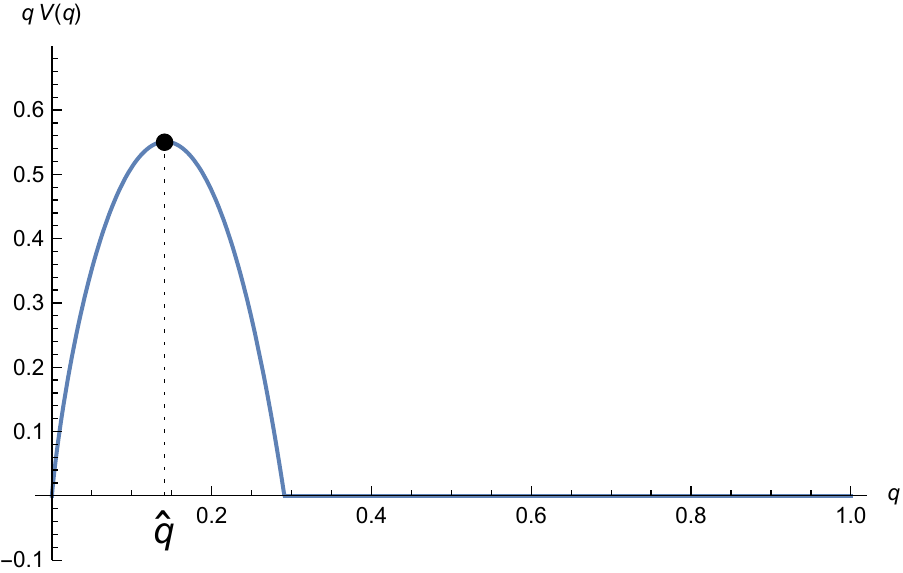}  
			\end{array}$
		\end{center}
	\caption{Graph of $q V(q)$ as a function of $q$. Parameters are $M=100, \, r=0.39, \sigma=4.1$ and $l=1.52$. The unique maximum is obtained at $\hq=0.14$. 
	}
	\label{qV_fig}
\end{figure}

Lastly, Figure \ref{qV_fig} graphically illustrates Proposition \ref{prop-critical-l}. As we check in the proof of the proposition, the map $q\mapsto q V(q) $ is strictly concave on $[0,p]$. The unique maximizer $\hat q$ in \eqref{hq_def} is marked in the figure.

\section{Concluding Remark}

In the cyber-security context, we develop a reputation game model between a suspect and a defender, and fully analyze the equilibrium interaction between them. As far as we know, our game model is the first to include the optimal termination of the game with asymmetric information, imperfect monitoring, and continuous-time Bayesian updates. Using the game model, we provide an empirical and theoretical methods of estimating the initial suspicion level.

As a future research, we plan to generalize our cyber-security game model by incorporating time-dependent noise size (periodic patterns of noise) and multidimensional signal processes (traffics from multiple channels).

\section*{Acknowledgement}
This work was supported by the National Research Foundation of Korea (NRF) grant funded by the Korea government (MSIT) (No. 2020R1C1C1A01014142, No. 2021R1I1A1A01050679, and No. 2021R1A4A1032924). 

\appendix

\addcontentsline{toc}{section}{Appendices}

\end{document}